\documentclass[aap,MSNbibl,citesort,dvips]{arximspdf}
\usepackage{graphicx}

\doi{10.1214/12-AAP842} 
\volume{22}
\issue{6}
\pubyear{2012}
\firstpage{2165}
\lastpage{2209}

\makeatletter
\newcommand{\IP}{\mathbb{P}}
\newcommand{\IE}{\mathbb{E}}
\newcommand{\IR}{\mathbb{R}}
\newcommand{\IZ}{\mathbb{Z}}
\newcommand{\IN}{\mathbb{N}}
\newcommand{\IT}{\mathbb{T}}

\newcommand{\rmP}{\mathrm{P}}
\newcommand{\rmE}{\mathrm{E}}
\newcommand{\rmT}{\mathrm{T}}

\newcommand{\tX}{\tilde{X}}

\newcommand{\hX}{\hat{X}}
\newcommand{\hS}{\hat{S}}
\newcommand{\hQ}{\hat{Q}}

\newcommand{\cO}{\mathcal{O}}
\newcommand{\cA}{\mathcal{A}}
\newcommand{\cL}{\mathcal{L}}

\newcommand\ind{\mathbf{1}}
\newcommand\e{\varepsilon}
\newcommand\integ[1]{\lfloor #1 \rfloor}

\newtheorem{thmm}{Theorem}[section]
\newtheorem{propn}[thmm]{Proposition}
\newtheorem{lemma}[thmm]{Lemma}
\newtheorem{lemmaa}{Lemma}
\newproclaim{eg}[thmm]{Example}
\newproclaim{defn}[thmm]{Definition}
\newproclaim{remark}[thmm]{Remark}
\newproclaim{notn}[thmm]{Notation}

\newproclaim{assumption}[thmm]{Assumption}
\makeatother

\begin{document}
\begin{frontmatter}

\title{The spatial $\Lambda$-Fleming--Viot process on a large torus:
Genealogies in the presence of recombination}
\runtitle{The spatial $\Lambda$-Fleming--Viot process with recombination}

\begin{aug}
\author[A]{\fnms{A. M.} \snm{Etheridge}\thanksref{t1}}
\and
\author[B]{\fnms{A.} \snm{V\'eber}\corref{}\thanksref{t2}\ead[label=e1]{amandine.veber@cmap.polytechnique.fr}}

\thankstext{t1}{Supported in part by EPSRC Grant EP/E065945/1.}
\thankstext{t2}{Supported in part by the \textit{chaire Mod\'elisation Math\'ematique et Biodiversit\'e} of Veolia Environnement-\'Ecole Polytechnique-Museum National d'Histoire Naturelle-Fondation X and by the ANR project MANEGE (ANR-09-BLAN-0215).}

\runauthor{A. M. Etheridge and A. V\'eber}

\affiliation{University of Oxford and CMAP---\'Ecole Polytechnique}

\address[A]{Department of Statistics\\
University of Oxford\\
1 South Parks Road\\
Oxford OX1 3TG\\
United Kingdom} 

\address[B]{CMAP---\'Ecole Polytechnique\\
Route de Saclay\\
91128 Palaiseau Cedex\\
France}
\end{aug}

\received{\smonth{6} \syear{2011}}
\revised{\smonth{1} \syear{2012}}

%
\begin{abstract}
We extend the spatial $\Lambda$-Fleming--Viot process introduced in
[\textit{Electron.~J. Probab.} \textbf{15} (2010) 162--216] to
incorporate recombination. The process models allele frequencies in a
population which is distributed over the two-dimensional torus $\IT(L)$
of sidelength $L$
and is subject to two kinds of reproduction events: \textit{small events}
of radius $\cO(1)$
and much rarer \textit{large events} of radius $\cO(L^{\alpha})$ for some
$\alpha\in(0,1]$.
We investigate the correlation between the times to the most recent
common ancestor of
alleles at two linked loci for a sample of size two from the population.
These individuals are initially sampled from ``far apart'' on the torus.
As $L$ tends to infinity, depending on the frequency of the large
events, the recombination rate
and the initial distance between the two individuals sampled, we obtain
either a complete
decorrelation of the coalescence times at the two loci, or a sharp
transition between a
first period of complete correlation and a subsequent period during
which the remaining times
needed to reach the most recent common ancestor at each locus are independent.
We use our computations to derive approximate probabilities of identity
by descent
as a function of the separation at which the two individuals are sampled.
\end{abstract}

%

\end{frontmatter}

\section{Introduction}\label{intro}

\subsection{Background}
In the 30 years since its introduction, Kingman's coalescent has become a
fundamental tool in population genetics. It provides an elegant\vadjust{\goodbreak}
description of the genealogical trees relating individuals in a
sample from a highly idealized biological population, in which it is assumed
that all individuals
are selectively neutral and experience identical conditions, and
that population size is constant.
Spurred on by the flood of DNA
sequence data, theoreticians have successfully extended the classical
coalescent to incorporate more realistic biological assumptions such as
varying population size, natural selection and genetic structure.
However, it has proved surprisingly difficult to produce satisfactory extensions
for populations living (as many do) in continuous two-dimensional
habitats---a problem dubbed \textit{the pain in the torus} by Felsenstein~\cite{Fel75}.

In the classical models of population genetics, it is customary to
assume that
populations are either panmictic, meaning, in particular, that they
have no
spatial structure, or that they are subdivided into ``demes.'' The demes sit
at the vertices of a graph which is chosen to caricature the geographic
region in which the population resides. Thus, for example, for a population
living in a two-dimensional spatial continuum one typically takes
the graph to be (a subset of) $\IZ^2$. Reproduction takes place within
demes and interaction between the subpopulations is through migration
along the edges of the graph. Models of this type are collectively known
as stepping stone models.

However, in order to apply a stepping stone model to
populations that are distributed
across continuous space, one is forced to make an artificial subdivision.
Moreover, the predictions of
stepping stone models
fail to match observed patterns of genetic variation. For example,
they overestimate genetic diversity (often by many orders of magnitude) and
they fail to predict the long-range correlations in allele frequencies seen
in real populations.

In recent work
\cite{etheridge2008,BEV2010,bartonkelleheretheridge2010}
we introduced a new framework in which to
model populations evolving in a spatial continuum. The key idea, which
enables us to overcome the pain in the torus, is that reproduction is
driven by a Poisson process of events which are based on geographical space
rather than on individuals. This leads, in particular,
to a
class of models that could reasonably be called \textit{continuum stepping stone
models}, but it also allows one to incorporate large-scale
extinction/recolonization events. Such events dominate the demographic
history of many species. They appear in our framework as
``local population bottlenecks.''
In~\cite{bartonkelleheretheridge2010}, we show
(numerically) how the
inclusion of such events can lead to long-range correlations in allele
frequencies. In~\cite{BEV2010} a rigorous mathematical analysis of a
class of
models on a torus in $\IR^2$ illustrates the reduction in genetic diversity
that can result from such large-scale
demographic events. We expand further on this in Section~\ref{subsibd}.
Thus, large-scale events
provide one plausible explanation of the two deficiencies of stepping
stone models highlighted above, but
of course they are not the only possible explanation.

A natural question now arises: how could we infer the existence of these
large-scale events from data?
One possible answer is through correlations
in patterns of variation at different genetic loci.
Recall that in a diploid population\vadjust{\goodbreak} (in which
chromosomes are carried in pairs)
correlations between linked genes (i.e., genes occurring on the
same chromosome) are broken down over time by recombination (which
results in two genes on the same chromosome being inherited from
different chromosomes in the parent). We say that genes are \textit{loosely}
linked, if the rate of recombination events is high [e.g., if
the chance of a recombination in a single
generation is ${\mathcal O}(1)$]. In the Kingman coalescent,
genealogies relating loosely linked genes evolve independently. This is
because on the timescale of the coalescent, the states in which
lineages ancestral to both loci are in the same individual vanish
instantaneously. It is well known that if a population experiences a
bottleneck, this is no longer the case. As we trace backward in time,
when we reach the bottleneck, we expect to see a significant proportion
of surviving lineages coalesce \textit{at the same time} and so we see
correlations in genealogies even at \textit{unlinked} loci.
With local bottlenecks we can expect
a rather more complicated picture. The degree of correlations across
loci will depend upon the spatial separation of individuals in the sample.

The purpose of this paper is to extend the model of~\cite{BEV2010} to
diploid populations, to incorporate recombination, and to
provide a first rigorous analysis of the
correlations in genealogies at different loci in the presence of
local extinction/recolonization events.
Since the questions we shall address
and some of the methods we shall use here are related to those of
\cite{BEV2010}, the reader may find it useful to have some familiarity with
the results of that paper.

\subsection{The model}
\label{subsmodel}
In~\cite{BEV2010}, we introduced the \textit{spatial $\Lambda$-Fleming--Viot process}
as a model of a haploid population evolving in a spatial continuum.
It is a Markov process taking
its values in the set of functions which associate to each point of the
geographical space a probability measure on a compact space, $K$, of genetic
types.
If $\Phi$ is the current state of the population and $x$ is
a spatial location, the measure $\Phi(x)$ can be interpreted as the
distribution of the type of an individual sampled from location $x$.
The dynamics are
driven by a Poisson point process of \textit{events}. An event specifies
a spatial region, $A$, say, and a number $u\in(0,1]$.
As a result of the event, a
proportion $u$ of individuals within $A$ are replaced
by offspring of a parent sampled from a point
picked uniformly at random from $A$. In~\cite{BEV2010} the regions $A$ are
chosen to be discs of random radius (whose centers fall with intensity
proportional to the Lebesgue measure) and the
distribution of $u$ can depend on the radius of the disc. Under
appropriate conditions,
existence and uniqueness in law of the process were established.

Here we wish to extend this framework in a number of directions.
First, whereas in~\cite{BEV2010} a single parent was chosen from the
region $A$,
here we allow $A$ to be repopulated by the offspring of a finite
(random) number
of its inhabitants. Second, we assume that the population is diploid.
We shall follow (neutral) genes at two
distinct (linked) loci, with recombination acting between them.
Writing $K_1$ and $K_2$ for the possible types at the two loci, the
type of an individual is an\vadjust{\goodbreak} element of $K_1\times K_2$
(which we can identify with $[0,1]\times[0,1]$).
As in~\cite{BEV2010}, we work on the torus $\IT(L)$ of side $L$ in
$\IR^2$
and we suppose that there are two types of events: small events, affecting
regions of radius $\cO(1)$, which might be thought of as ``ordinary''
reproduction events; and ``large'' events, representing extinction/recolonization
events, affecting regions of radius
$\cO(L^\alpha)$ where $\alpha\in(0,1]$ is a fixed parameter.
In order to
keep the notation as simple as possible, we shall only allow two different
radii for our events, $R_s$ corresponding to ``small'' reproduction
events and
$R_B L^\alpha$ corresponding to ``large'' local bottlenecks.
We shall also suppose
that the corresponding
proportions $u_s$ and $u_B$ are fixed. Neither of these assumptions is
essential to the results, which would carry over to the more general setting
in which each of $R_s$, $R_B$, $u_s$ and $u_B$ is sampled (independently)
from given distributions each time an event occurs.

Let us specify the dynamics of the process more precisely.
Let:
\begin{itemize}
\item$R_s,R_B\in(0,\infty)$, $u_s,u_B \in(0,1)$ and $\alpha\in(0,1]$;
\item$\lambda_s,\lambda_B$ be two distributions on $\IN=\{
1,2,\ldots\}$
with bounded support and such that $\lambda_s(\{1\})<1$;
\item$(\rho_L)_{L\in\IN}$ be an increasing sequence such that
$\rho_L\geq\log L$ for all $L\in\IN$, and $L^{-2\alpha}\rho_L$
tends to a finite limit (possibly zero) as $L\rightarrow\infty$; and
\item$(r_L)_{L\in\IN}$ be a nonincreasing sequence with values in $(0,1]$.
\end{itemize}
For $L\in\IN$, we denote by
$\Pi_s^L$ a Poisson point process on $\IR\times\IT(L)$ with intensity
measure $dt\otimes dx$, and by $\Pi_B^L$ another Poisson point process on
$\IR\times\IT(L)$, independent of $\Pi^L_s$, with intensity measure
$(\rho_LL^{2\alpha})^{-1}dt\otimes dx$.
The spatial $\Lambda$-Fleming--Viot process
$\Phi^L$ on $\IT(L)$ evolves as follows.

\textit{Small events}: If $(t,x)$ is a point of $\Pi^L_s$, a
reproduction event
takes place at time~$t$ within the closed ball $B(x,R_s)$:
\begin{itemize}
\item a number $j$ is sampled according to the measure $\lambda_s$;
\item $j$ sites, $z_1,\ldots,z_j$ are selected uniformly at random from
$B(x,R_s)$; and,
\item for each $i=1,\ldots,j$, a type $(a_i,b_i)$ is sampled
according to $\Phi_{t-}^L(z_i)$.
\end{itemize}
If $j>1$, then for all $y\in B(x,R_s)$,
\[
\Phi_t^L(y):= (1-u_s)\Phi_{t-}^L(y) +\frac{u_s(1-r_L)}{j}\sum
_{i=1}^{j}\delta_{(a_i,b_i)} +\frac{u_sr_L}{j(j-1)}\sum_{i_1\neq
i_2}\delta_{(a_{i_1},b_{i_2})}.
\]
If $j=1$,
for each $y\in B(x,R_s)$,
\[
\Phi_t^L(y):= (1-u_s)\Phi_{t-}^L(y) +u_s\delta_{(a_1,b_1)}.
\]
In both cases, sites outside $B(x,R_s)$ are not affected.

\textit{Large events}: If $(t,x)$ is a point of $\Pi^L_B$, an
extinction/recolonization
event takes place at time $t$ within the closed ball $B(x,L^{\alpha}R_B)$:
\begin{itemize}
\item a number $j$ is sampled according to the measure $\lambda_B$;
\item $j$ sites, $z_1,\ldots,z_j$ are selected uniformly at random from
$B(x,L^{\alpha}R_B)$; and
\item for each $i=1,\ldots,j$, a type $(a_i,b_i)$ is sampled
according to $\Phi_{t-}^L(z_i)$.\vadjust{\goodbreak}
\end{itemize}
For each $y \in B(x,L^{\alpha}R_B)$,
\[
\Phi_t^L(y):= (1-u_B)\Phi_{t-}^L(y) +\frac{u_B}{j} \sum
_{k=1}^j\delta
_{(a_k,b_k)}.
\]
Again, sites outside the ball are not affected.
%
\begin{remark}
\label{remarksonmodel}
(1) The scheme of choosing $j$ parental locations and then
sampling a parental type at each of those locations is convenient when
one is interested in tracing lineages ancestral to a sample from the
population.
It can be thought of as sampling $j$ individuals, uniformly at random
from the ball
affected by the reproduction (or extinction/recolonization) event, to
reproduce. Of course this scheme allows for the possibility of more
than one parent contributing offspring so that we should more
correctly call this model a spatial $\Xi$-Fleming--Viot process, but to
emphasize the
close link with previous work we shall abuse terminology and use the name
$\Lambda$-Fleming--Viot process.\vspace*{-6pt}
\begin{longlist}[(2)]
\item[(2)] The recombination scheme mirrors that generally employed in Moran
models. The quantity $r_L$ is the proportion of offspring who, as a
result of recombination, inherit the types at the two loci from different
parental chromosomes. We have chosen to
sample the types of those two chromosomes
from different points in space. The result of this is that provided
the individuals sampled from the current population are in distinct
geographic locations,
if two ancestral lineages
are at spatial distance zero, then they are necessarily in the same individual.
This is mathematically convenient (cf. Remark~\ref{remmarkovproperty}) but,
arguably, not terribly natural biologically. However,
changing the sampling scheme, for example, so that the two recombining
chromosomes are sampled from the same location,
would not materially change our results.

\item[(3)]
We are assuming that recolonization is so rapid after an extinction
event that the effects of recombination during recolonization are
negligible.
\item[(4)]
Since $\IT(L)$ is compact, the overall rate at which events fall is finite
for any~$L$ and the corresponding \textit{spatial $\Lambda$-Fleming--Viot process
with recombination} is well-defined. Notice that a given site, $x$, say,
is affected
by a small event at rate $\pi R_s^2 = \cO(1)$ (since the center of the event
must fall within a distance $R_s$ of~$x$),
whereas it is hit by a large event at rate
$\pi R_B^2\rho_L^{-1} = \cO(\rho_L^{-1})$. So reproduction events are
frequent, but massive extinction/recolonization events are rare.
\end{longlist}
\end{remark}

\subsection{Genealogical relationships}
\label{genealogies}
Having established the (forward in time) dynamics of allele frequencies
in our model, we now turn to the genealogical relationships between
individuals in a sample from the population.

First suppose that we are tracing the lineage ancestral to a single
locus on
a chromosome carried by just one individual in the current population.
Recombination does\vadjust{\goodbreak} not affect us and we see that the lineage will move
in a series of jumps: if its current location is $z$, then it will jump
to $z+x$ (resp., $z+L^\alpha x$) due to a small (resp., large) event with
respective intensities
%
%
\begin{equation}\label{intensities}
L_{R_s}(0,x)  u_s  \frac{dx}{\pi R_s^2}\quad \mbox{and}\quad
\frac
{L_{R_B}(0,x)}{\rho_L}  u_B  \frac{dx}{\pi R_B^2},
\end{equation}
where $L_R(x,y)$ denotes the volume of the intersection
$B(x,R)\cap B(y,R)$ [viewed as a subset of $\IT(L)$ for the first
intensity measure, and of $\IT(L^{1-\alpha})$ for the second].
To see this, note first that by translation invariance of the model we may
suppose that $z=0$. In order for the lineage to experience a small
jump, say,
from the origin to $x$, the origin and the position $x$ must be
covered by the same event. This means that the center of the event must lie
in both $B(0,R_s)$ and $B(x,R_s)$. The rate at which such events occur
is $L_{R_s}(0,x)$. The lineage will only jump if it is sampled from the
portion $u_s$ of the population that are offspring of the event and
then it will jump to the position of its parent, which is uniformly
distributed on a ball of area $\pi R_s^2$. Combining these observations gives
the first intensity in~(\ref{intensities}).
A lineage ancestral to a single locus in a single individual
thus follows a compound Poisson process on $\IT(L)$.

Suppose now that we sample a single individual, but trace back its ancestry
at \textit{both} loci. We start with a single lineage which moves, as
above, in
a series of jumps as long as it is in the fraction $u_s(1-r_L)$ of
``nonrecombinants'' in the population. However, every time it is hit by a
small event, there is a probability $u_sr_L$ that it was created by
recombination from two parental chromosomes, whose locations are sampled
uniformly at random from the region affected by the event.
If this happens, we must follow two distinct lineages, one for each locus,
which jump around $\IT(L)$ in an a priori correlated manner
(since they may be hit by the same events),
until they coalesce again. This will happen if they are both affected
by an
event (small or large) and are both derived from the same parent (which
for a given
event
has probability $1/j$ in our notation above).

Thus, the ancestry of the two loci from our sampled individual
is encoded in a system of splitting and coalescing lineages. If we now
sample two individuals, $(A,B)$ and $(a,b)$, we represent their
genealogical relations at the two loci by a process $\cA^L$ taking
values in the set of partitions of $\{A,a,B,b\}$ whose blocks are
labeled by an element of $\IT(L)$. As in~\cite{BEV2010}, at time
$t\geq
0$ each block of $\cA^L_t$ contains the labels of all the lineages
having a common ancestor (i.e., carried by the same individual) $t$
units of time in the past, and the mark of the block records the
spatial location of this ancestor. The only difference with the
ancestral process defined in~\cite{BEV2010} is that blocks can now
split due to a recombination event.

Of course, if $r_L$ is small, then the periods of time when the lineages
are in a single individual, that is, during which
they have coalesced and not split apart
again, can be rather extensive. This has the potential to create
strong correlations between the two loci. The other source of correlation
is the large events which can cause coalescences between lineages even when
they are geographically far apart.
To gain an understanding of these correlations, we ask the following question:

\textit{The problem}:
Given $\alpha$, $\rho_L$ and $r_L$, is there a minimal distance
$D_L^*$ such that, asymptotically as $L\rightarrow\infty$,
\begin{itemize}
\item if we sample two individuals $(A,B)$ and $(a,b)$ at distance at
least $D_L^*$ from each other, then the coalescence time of the
ancestral lineages of $A$ and $a$ is independent of that of the ancestral
lineages of $B$ and $b$ (in other words, genealogies at the two loci are
completely decorrelated);
\item if two individuals are sampled at a distance less than $D_L^*$,
then the genealogies at the two loci are correlated (i.e., the lineages
ancestral to $A$ and $B$, and, similarly, those of $a$ and $b$,
remain sufficiently ``close together'' for a sufficiently long time that
there is a significant chance that
the coalescence
of $A$ and $a$ implies that of~$B$ and $b$ at the same time or soon after)?
\end{itemize}

\subsection{Main results}
\label{results}
Before stating our main results, we introduce some notation.
We shall always denote the types of the two individuals in our
sample by $(A,B)$ and $(a,b)$.
The same letters will be used to distinguish the corresponding ancestral
lineages.
As we briefly mentioned in the last section, the genealogical
relationships between the two loci at time $t\geq0$ before the
present are represented by a marked partition of $\{A,a,B,b\}$, in which
each block corresponds to an individual in the ancestral population at
time $t$
who carries lineages ancestral to our sample.
The labels in the block are those of the corresponding lineages and the mark
is the spatial location of the ancestor. For any such marked partition
$a_L$, we write $\IP_{a_L}$ for the probability measure under which the
genealogical process starts from $a_L$, with
the understanding that marks then evolve on the torus $\IT(L)$.
Typically, our initial configuration will be of the form
\[
a_L:= \{(\{A,B\},x_L^1),(\{a,b\},x_L^2)\},
\]
where the separation $x_L:= x_L^1-x_L^2$ between the two sampled
individuals will be assumed to be large.
The coalescence times of the ancestral lineages at each locus are
denoted by $\tau_{Aa}^L$ and $\tau_{Bb}^L$. Finally, we write $|x|$
for the
Euclidean norm of $x\in\IR^2$ (or in a torus of any size) and $\sigma
^2>0$ is a constant, whose value is given just after~(\ref{defsigma}).
(It corresponds, after a suitable space--time rescaling, to the limit as
$L\rightarrow\infty$ of the variance of the displacement of a
lineage during a time interval of length one.)

For later comparison, we first record the
asymptotic behaviour of the coalescence time at a single locus. The
proof of the following result is in Section~\ref{section1locus}.
%
\begin{propn}\label{prop1locus}
Suppose that for each $L\in\IN$ the two individuals comprising
our initial configuration $a_L$ are
at separation $x_L\in\IT(L)$.
Suppose also that $\frac{\log|x_L|}{\log L}\rightarrow\beta\in
(\alpha,1]$
as $L\rightarrow\infty$. (In particular, $\alpha<1$ here.)
Then:
\begin{longlist}
\item[(a)] For all $t\in[\beta,1]$,
\[
\lim_{L\rightarrow\infty} \IP_{a_L}\bigl[\tau_{Aa}^L> \rho
_LL^{2(t-\alpha)}\bigr] = \frac{\beta-\alpha}{t-\alpha}.
\]
\item[(b)] For all $t>0$,
\[
\lim_{L\rightarrow\infty} \IP_{a_L}\biggl[\tau_{Aa}^L >\frac
{1-\alpha
}{2\pi\sigma^2} \rho_L L^{2(1-\alpha)}\log L t \biggr] = \frac
{\beta
-\alpha}{1-\alpha}  e^{-t}.
\]
\end{longlist}
\end{propn}

\begin{remark}Observe that the timescale considered in case (b) above
coincides with the quantity $\varpi_L$ defined in Theorem~3.3 of \cite
{BEV2010}. Indeed, using the notation of~\cite{BEV2010}, the variance
$\sigma^2$ is given by the following limit:
\[
\sigma^2=\lim_{L\rightarrow\infty} \frac{\rho_L}{L^{2\alpha}}
\sigma
_s^2 + \sigma_B^2,
\]
where $\sigma_s^2$ and $\sigma_B^2$ are defined in equation~(20) of
\cite{BEV2010}. Now, if $\rho_LL^{-2\alpha}\rightarrow0$ as in case
(a) of Theorem~3.3, we obtain
\[
\frac{1-\alpha}{2\pi\sigma^2} \rho_L L^{2(1-\alpha)}\log L
\approx
\frac{1-\alpha}{2\pi\sigma_B^2} \rho_L L^{2(1-\alpha)}\log L,
\]
while if $\rho_LL^{-2\alpha} \rightarrow1/b>0$, we have
\begin{eqnarray*}
\frac{1-\alpha}{2\pi\sigma^2} \rho_L L^{2(1-\alpha)}\log L
&\approx&
\frac{1-\alpha}{2\pi(({1}/{b}) \sigma^2_s +\sigma_B^2)}
\frac
{1}{b}  L^2\log L \\
&=& \frac{1-\alpha}{2\pi(\sigma^2_s +b\sigma
_B^2)}
L^2\log L .
\end{eqnarray*}
In both cases, the timescale considered in Proposition~\ref{prop1locus} is the same
as the quantity $\varpi_L$ of Theorem~3.3 of~\cite{BEV2010}.
\end{remark}

In the case $\alpha=0$, Proposition~\ref{prop1locus} precisely
matches corresponding results of~\cite{CG1986}
and~\cite{ZCD2005} for coalescing random walks on a torus in $\IZ^2$.
For $\alpha>0$, we
see that if lineages start at a separation of $\cO(L^\beta)$, with
$\beta>\alpha$, then the small events don't affect the asymptotic
coalescence times; they are the same as those for
a random walk with bounded jumps on $\IT(L^{1-\alpha})$ started at
separation $\cO(L^{\beta(1-\alpha)})$.
In particular, the first statement tells us that the chance that
coalescence occurs at a time
$\ll\rho_L L^{2(1-\alpha)}\log L$ is $(1-\beta)/(1-\alpha)$. If this
does not happen, then since the time taken
for the random walks to reach their equilibrium distribution is
$\cO(\rho_L L^{2(1-\alpha)}\log L)$, in these units, the additional
time that we must
wait to see a coalescence is asymptotically exponential.

When we consider the genealogies at two loci, several regimes appear
depending on the recombination rate and the initial distance between
the individuals sampled.
%
\begin{thmm}\label{theocorrelation}
Suppose $(a_L)_{L\in\IN}$ is as in Proposition~\ref{prop1locus}. If
%
%
\begin{equation}\label{cond}
\limsup_{L\rightarrow\infty} \frac{\log(1+{\log\rho
_L}/{(r_L\rho_L)})}{2\log L}\leq\beta-\alpha,
\end{equation}
then we have the following:

\begin{longlist}
\item[(a)] For all $t\in[\beta,1]$,
\[
\lim_{L\rightarrow\infty} \IP_{a_L}\bigl[\tau_{Aa}^L \wedge\tau
_{Bb}^L>\rho_L L^{2(t-\alpha)}\bigr] = \frac{(\beta-\alpha
)^2}{(t-\alpha)^2}.
\]
\item[(b)] For all $t>0$,
\[
\lim_{L\rightarrow\infty} \IP_{a_L}\biggl[\tau_{Aa}^L \wedge\tau
_{Bb}^L >\frac{1-\alpha}{2\pi\sigma^2} \rho_L L^{2(1-\alpha
)}\log L t\biggr]
= \frac{(\beta-\alpha)^2}{(1-\alpha)^2}  e^{-2t}.
\]
\end{longlist}
\end{thmm}

Under the conditions of Theorem~\ref{theocorrelation}, the
individuals are initially sampled at a distance much larger than the
radius of the large events, and recombination is fast enough for the
coalescence times at the two loci to be asymptotically independent (see
Remark~\ref{remindependence}).
For slower recombination rates this is no longer the case.
When Condition~(\ref{cond}) is not satisfied, we have instead:
%
\begin{thmm}\label{theosemi-correl}
Suppose $(a_L)_{L\in\IN}$ is as in Proposition~\ref{prop1locus}.
Assume there exists $\gamma\in(\beta,1)$ such that
%
%
\begin{equation}\label{cond2}
\lim_{L\rightarrow\infty}\frac{\log(1+{\log\rho
_L}/{(r_L\rho
_L)})}{2\log L}= \gamma-\alpha.
\end{equation}
Then:
\begin{longlist}
\item[(a)] For all $t\in[\beta,\gamma]$,
\[
\lim_{L\rightarrow\infty} \IP_{a_L}\bigl[\tau_{Aa}^L\wedge\tau
_{Bb}^L>\rho_L L^{2(t-\alpha)}\bigr] = \frac{\beta-\alpha
}{t-\alpha}.
\]
\item[(b)] For all $t\in(\gamma,1]$,
\[
\lim_{L\rightarrow\infty} \IP_{a_L}\bigl[\tau_{Aa}^L\wedge\tau
_{Bb}^L>\rho_L L^{2(t-\alpha)}\bigr]= \frac{(\beta-\alpha)(\gamma
-\alpha
)^2}{(\gamma-\alpha)(t-\alpha)^2}.
\]
\item[(c)] For all $t>0$,
\[
\lim_{L\rightarrow\infty} \IP_{a_L}\biggl[\tau_{Aa}^L\wedge\tau_{Bb}^L
>\frac{1-\alpha}{2\pi\sigma^2} \rho_L L^{2(1-\alpha)}\log L t\biggr]
= \frac{(\beta-\alpha)(\gamma-\alpha)^2}{(\gamma-\alpha
)(1-\alpha)^2} e^{-2t}.
\]
\end{longlist}
\end{thmm}

This time, we observe a ``phase transition'' at time $\rho_LL^{2(\gamma
-\alpha)}$.
Asymptotically, coalescence times are completely correlated for times of
$\cO(\rho_LL^{2(\gamma-\alpha)})$, but conditional on being greater
than this ``decorrelation threshold'' they are independent. To understand
this threshold, recall from
Proposition~\ref{prop1locus} that, initially, coalescence of
lineages ancestral to a single locus happens on the exponential
timescale $\rho_LL^{2(t-\alpha)}$, $t\in[\beta,1]$
and is driven by large events.
This tells us that the effect of recombination will be felt only if
exactly one of the lineages ancestral to $A$ and~$B$ (or to $a$ and
$b$) is ``hit'' by
a large event. Since recombination events between $A$ and $B$ result in
only a small separation of the corresponding ancestral lineages,
we can expect that many of them will rapidly be followed by coalescence
of the corresponding
lineages (due to small events). This leads us to the idea of an
``effective'' recombination event, which is one following which at least
one of the
lineages ancestral to $A$ and $B$ is affected by a large event {\em
before} they coalesce due to small events. We shall see in Proposition
\ref{propeffectiverecomb} that
recombination is ``effective'' on the linear timescale $\rho_L(1+(\log
\rho_L)/(r_L\rho_L))  t$, $t\geq0$.
Under condition~(\ref{cond2})
the timescales of coalescence and effective recombination cross over
precisely at time $\rho_LL^{2(\gamma-\alpha)}$.

Two cases remain:

\textit{The case $\alpha<\beta\leq1$, $\gamma\geq1$}: If $\gamma>1$,
the arguments of the proof of Theorem~\ref{theosemi-correl} show
that the
recombination is too slow to be effective on the timescale of
coalescence and so the coalescence times at the two loci are completely
correlated and are given by
Proposition~\ref{prop1locus}.
For $\gamma=1$, the result depends on the precise form of $(\log\rho
_L)/(r_L\rho_L)$.
If it remains close enough to $L^{2(1-\alpha)}$ (or smaller), the proof
of Theorem~\ref{theosemi-correl} shows
that lineages are completely correlated on the timescale $\rho
_LL^{2(t-\alpha)}$, $t\leq1$, and then, conditional on not having
coalesced before $\rho_LL^{2(1-\alpha)}$, they evolve independently on
the timescale $\rho_LL^{2(1-\alpha)}\log L~t$.
On the other hand, if $(\log\rho_L)/(r_L\rho_LL^{2(1-\alpha)})$ grows
to infinity sufficiently fast, then, just as in the case $\gamma>1$,
recombination is too slow to be effective.

\textit{The case $\beta\leq\alpha\leq1$}: If we drop our assumption
that the separation of the individuals in our sample is much
greater than the radius of the largest events, then we can no longer
make such precise statements. Proposition~6.4(a) in~\cite{BEV2010} (with
$\psi_L=L^{\alpha}$) shows that the coalescence time for lineages
ancestral to a single locus will now be at most $\cO(\rho_L)$.
This does tell us that if $r_L\rho_L\rightarrow0$ as $L\rightarrow
\infty$, then
asymptotically we will not see any
recombination before coalescence and the coalescence times $\tau
_{Aa}^L$ and $\tau_{Bb}^L$ are identical. However, in contrast
to the setting of Proposition~\ref{prop1locus}, even asymptotically,
their common value depends on the exact
separation of the individuals sampled.
The same reasoning is valid when $r_L\ll(\log\rho_L)/\rho_L$. In this
case, although
we may see some recombination events before any coalescence occurs, a
closer look at the proof of Proposition~\ref{propeffectiverecomb} reveals
that the time spent in distinct individuals by the two lineages
ancestral to $A,B$, say,
in $\cO(\rho_L)$ units of time, is negligible
compared to $\rho_L$. Thus, with high probability, any large event
affecting lineages ancestral
to our sample will occur at a time when the lineages ancestral to $A$
and $B$ are in the same
individual, (as are those
ancestral to $a$ and $b$).
As a result, once again $\tau_{Aa}^L=\tau_{Bb}^L$ with probability
tending to $1$.

On the other hand, suppose $r_L$ remains large enough that lineages
ancestral to $A$ and $B$
have a chance to be hit by a large event while they are in different
individuals and
thus jump to a separation $\cO(L^{\alpha})$ (the \textit{effective
recombination} of
Section~\ref{subseffectiverecomb}). We are still unable to recover
precise results. The reason is
that even after such an event, we may be in a situation in which all
lineages could be hit by the same large event, or at
least remain at separations $\cO(L^{\alpha})$. But we shall see that a
key to the proofs of
Theorems~\ref{theocorrelation} and~\ref{theosemi-correl} is the
fact that, in the settings
considered there, where individuals are sampled from far apart,
whenever two
lineages come to within $2R_BL^{\alpha}$ of one another, the other
ancestral lineages
are still very far from them. This gives the pair time to merge without
``interference'' from
the other lineages. Since lineages at separations $\cO(L^\alpha)$ are
correlated and
their coalescence times depend strongly on their
precise (geographical) paths on this scale, it is difficult to quantify
the extent
to which the fact that the ancestral lines of $A,B$ and of $a,b$ start
within the same
individuals makes the coalescence times $\tau_{Aa}^L$ and $\tau_{Bb}^L$
\textit{more}
correlated. Nonetheless, this is an important question and will be
addressed elsewhere.

To answer our initial question, we see from Theorems~\ref{theocorrelation}
and~\ref{theosemi-correl} (and the subsequent discussion) that
$D_L^*$ is informally given by
\[
\log\biggl(1+\frac{\log\rho_L}{r_L\rho_L}\biggr)\approx2(\log
D_L^*-\alpha\log L),\qquad \mbox{i.e., }  D_L^*\approx L^\alpha
\sqrt
{1+\frac{\log\rho_L}{r_L\rho_L}}.
\]
When the sampling distance is greater than the radius of the largest
events, correlated
genealogies are only possible when recombination is slow enough, or
large events occur rarely enough,
that $(\log\rho_L)/(r_L\rho_L)\gg1$. If, for instance, $r_L\equiv
r>0$, the two loci are
always asymptotically decorrelated. On the other hand, if $\gamma$ is
as in~(\ref{cond2})
[note that $\gamma$ does not need to exist for condition (\ref{cond})
to hold] and the
sampling distance is $L^\beta$, Theorem~\ref{theocorrelation} shows
that if
$\beta\geq\gamma$, the genealogies at the two loci are asymptotically
independent,
whereas Theorem~\ref{theosemi-correl} tells us that if $\beta\in
(\alpha,\gamma)$,
there is a first phase of complete correlation. Thus,
$D_L^*\mbox{``}=\mbox{''}L^\gamma$.

Before closing this section, let us make two remarks:
%
\begin{remark}[(Bounds on the rates of large
events)]\label{remrates}
Recall that we imposed the condition $\log L\leq\rho_L\leq
CL^{2\alpha
}$. The reason for
the upper bound is that in~\cite{BEV2010}, we showed that the
coalescence of the ancestral lineages
is then driven by the large events and, moreover, is very rapid once
lineages are at separation
$\cO(L^{\alpha})$ (see the proof of Theorem~3.3 in~\cite{BEV2010}).
Similar results should hold,
although on different timescales,
in the other cases presented in~\cite{BEV2010}. However, to keep the
presentation of our results
as simple as possible, we have chosen to concentrate on this upper bound.
The (rather undemanding) lower bound is needed in the proof of
Proposition~\ref{propdecorr}.
\end{remark}

\begin{remark}[(Generalization of
Theorems~\protect\ref{theocorrelation} and~\protect\ref{theosemi-correl} to
distinct coalescence times)]\label{remindependence}
In these two theorems, we could also consider the probabilities of
events of the
form $\{\tau_{Aa}^L>\rho_LL^{2(t-\alpha)} \mbox{ and } \tau
_{Bb}^L>\rho_LL^{2(t'-\alpha)}\}$,
with $t<t'$. However, they can be computed by a simple application of
Theorem~\ref{theocorrelation} or~\ref{theosemi-correl} at time $t$,
and the Markov property.
Indeed, arguments similar to those of the proofs of Lemma~\ref{lemC} and
Lemma~\ref{lemmasep} in
Section~\ref{section1locus} tell us that the distance between
lineages ancestral to~$B$ and $b$ at time $\rho_LL^{2(t-\alpha)}$, conditional on not having
coalesced by this
time, lies in $[L^t/(\log L),L^t\log L]$. Proposition~\ref{prop1locus} then enables us to
conclude. We leave this generalization to the reader.
\end{remark}

The rest of the paper is laid out as follows. In Section~\ref{subsibd} we provide more detail
of the
motivation for the question addressed here. In Section~\ref{section1locus} we prove
Proposition~\ref{prop1locus} and collect several results on
genealogies of a sample from a
single locus that we shall need in the sequel. Since most of these
results are close
to those established in~\cite{BEV2010}, or require techniques used in
\cite{CG1986} and
\cite{ZCD2005} for similar questions on the discrete torus, their
proofs will only be sketched.
Our main results are proved in Section~\ref{section2loci}: we define
an effective
recombination rate in Section~\ref{subseffectiverecomb}, use it
to find an upper bound on the time we must wait before the two lineages
ancestral to $A$ and $B$
start to evolve independently
in Section~\ref{subsdecorrelation}
and finally derive the asymptotic coalescence times of our two pairs of
lineages in
Section~\ref{subsproofs}.

\section{Biological motivation}\label{subsibd}

In this section we expand on the biological motivation for our work.

It has long been understood that for many models of spatially
distributed populations, if
individuals are sampled sufficiently far from one another, then the
genealogical tree that records
the relationships between the alleles carried by those individuals at a
single locus
is well-approximated by a Kingman coalescent with
an ``effective population size'' capturing the influence of the geographical
structure. If the underlying population model is a stepping stone
model, with
the population residing in discrete demes located at the vertices of
$\IZ^2$ or $\IT(L)\cap\IZ^2$, individuals reproducing within demes
and migration modeled as a random walk, then the genealogical trees relating
individuals in a finite sample from the population are traced out by a
system of coalescing random walks.
The case in which random walks coalesce instantly on meeting
corresponds (loosely) to
a single individual living in each deme in which case the stepping
stone model
reduces to the voter model.
In this setting, and with symmetric nearest neighbour migration,
convergence to the Kingman coalescent as the separation of individuals
in the initial sample
tends to infinity
was established for $\IZ^2$ in~\cite{CG1986,CG1990}, and for $\IT
(L)\cap\IZ^2$ in~\cite{COX1989}.
In~\cite{CD2002,ZCD2005}, Z\"ahle, Cox and Durrett prove the same kind
of convergence for coalescing random walks on
$\IT(L)\cap\IZ^2$ with finite variance jumps and delayed coalescence
(describing the genealogy for a sample from Kimura's
stepping stone model on the discrete torus in which reproduction within
each deme is modeled by a Wright--Fisher diffusion).
In~\cite{limicsturm2006}, Limic and Sturm prove the analogous result
when mergers between random walks within a deme
are not necessarily pairwise.
In the same spirit but on the continuous space $\IT(L)$ and with
additional \textit{large} extinction/recolonization events
(similar to those described in Section~\ref{subsmodel}), the same
asymptotic behavior is obtained in~\cite{BEV2010} for the systems
of coalescing compound Poisson processes describing the genealogy of a
sample from
the spatial $\Lambda$-Fleming--Viot process, under suitable conditions
on the frequency and extent of the large events.

In all of these examples, the result stems from a separation of timescales.
For example, in~\cite{BEV2010} we were concerned with the genealogy of
a sample picked
uniformly at random from the whole torus.
Under this assumption, the time that two lineages need to be ``gathered''
close enough
together that they can both be affected by the same event
dominates the additional time the lineages take to coalesce, having
being gathered.
As explained in Section~\ref{results}, this decomposition does not hold
when lineages start
too close together, and so the tools
developed for well-separated samples are of no use in the study of
local correlations.
However, although we still cannot make precise statements about the genealogy
of samples which are initially too close together,
the work of Sections~\ref{subseffectiverecomb} and~\ref{subsdecorrelation}, which
are concerned with ``effective recombination'' and ``decorrelation,''
provides a
much better understanding than we had before of the local mechanisms that
create correlations between nearby lineages, how strong these
correlations are,
and how to ``escape'' them.

Our main results in this paper are concerned with samples taken at
``intermediate'' scales.
Individuals are sampled at pairwise distances much larger than the
radius of the largest
events, but these distances can still be much less than the radius of
the torus.
In this case, the ``gathering time'' of two lineages starting at
separation $x_L$
depends on that separation, but asymptotically this dependence is only
through $\log|x_L|/\log L$. As in the case of a uniform sample,
the gathering time dominates the additional
time to coalescence.
In Theorem~3.3 of~\cite{BEV2010} we showed that if we sample a finite number
of individuals uniformly at random from
the geographic range of a population which is subject to small and large
demographic events, then measuring time in units of size
$\varpi_L=\frac{1-\alpha}{2\pi\sigma^2}(\rho_L/L^{2\alpha
})L^2\log L$
(under the assumption on $\rho_L$ used here), their genealogical tree
is determined by Kingman's coalescent.
In particular, if $\rho_L<L^{2\alpha}$ (i.e., large events are not too
rare), one major effect of the presence of large
extinction/recolonization events is to reduce the \textit{effective
population size} and, consequently, genetic diversity.
The assumption of uniform sampling guarantees that initially ancestral
lineages are $\cO(L/\log L)$ apart.
Proposition~\ref{prop1locus} extends the result by showing that, if
we sample our individuals from much closer together, then we
should consider two timescales. The first is $(\rho_L/L^{2\alpha
})L^{2t},  t\in[\beta,1]$. The second kicks in after $\cO(\rho_L
L^{2(1-\alpha)})$, when
the lineages start to feel the fact that space is limited and their
ancestries evolve on the
linear timescale $\varpi_Lt$.
Now, by the same reasoning, if there were no large events, these
timescales would be, respectively,
$L^{2t},  t\in[\beta,1]$, and $\frac{1}{2\pi\sigma_s^2}L^2\log L
t$, $t>0$.
Of course, one never observes genealogies directly and so, for illustration,
we introduce (infinitely many alleles) mutation
into our model and compute the probability that two individuals sampled
at a given separation are \textit{identical by descent} (IBD) as a
function of the
exponent $\beta$.
In other words, what is the probability that the two individuals carry
the same type (at a given locus) because it was inherited from a common
ancestor.

\begin{figure}

\includegraphics{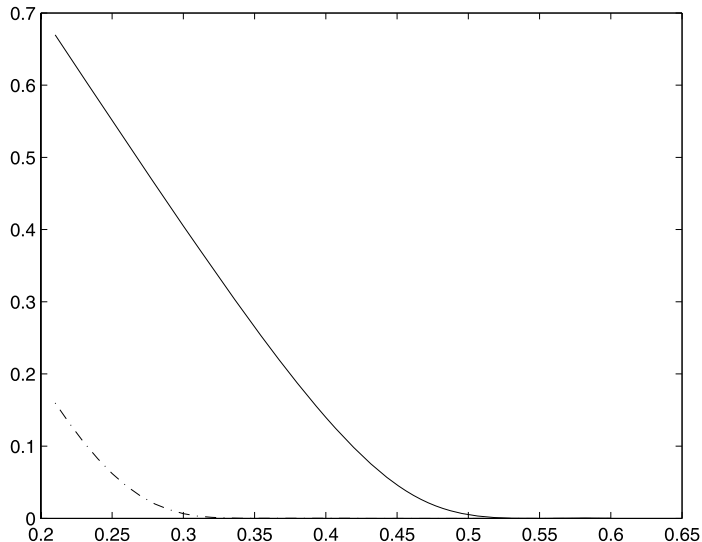}

\caption{Probability of IBD at a single locus, as a function of $\beta$.
Here, $L=10^5$, $\alpha=0.1$, $c_L=0.01$ and $\theta=10^{-3}$. The
solid line corresponds to the case with small and large events, the
dash--dot line to the case with only small events. Geographical
correlations vanish around $\beta=0.32$ without large events, and are
positive up to $\beta=0.52$ when large events occur.}
\label{fig1locus}
\end{figure}

Since mutations are generally assumed to occur at a linear rate, while
the first phase of the genealogical tree
develops on a much slower exponential timescale, for a given time
parameter $t\in[\beta,1]$, asymptotically as $L\rightarrow\infty$,
we would see either zero or infinitely many mutations on the tree.
However, let us suppose that~$L$ is large and write $\theta$ for the
mutation rate at locus $A$. We denote by $c_L$ the ratio $\rho
_L/L^{2\alpha}$.
Since IBD is equivalent to our individuals experiencing no mutation
between the time of their most recent common ancestor and the present,
the probability of IBD of two individuals sampled at distance $L^\beta$
is given by
\begin{eqnarray}\label{eqibd1}
\IE_{L^\beta}[e^{-2\theta\tau_{Aa}^L}] & \approx&\IE
_{L^\beta
}\bigl[e^{-2\theta\tau_{Aa}^L}\ind_{\{c_LL^{2\beta}\leq\tau
_{Aa}^L\leq
c_LL^2\}}\bigr] + \IE_{L^\beta}\bigl[e^{-2\theta\tau_{Aa}^L}\ind
_{\{
c_LL^2 < \tau_{Aa}^L\}}\bigr]\nonumber\\
& = &\int_{c_LL^{2\beta}}^{c_LL^2}e^{-2\theta t}\IP_{L^\beta}
[\tau
_{Aa}^L\in dt] + \int_{c_LL^2}^{\infty}e^{-2\theta t}\IP
_{L^\beta
}[\tau_{Aa}^L\in dt]\\
& \approx&(\beta-\alpha)\int_\beta^1 \frac{e^{-2\theta
c_LL^{2u}}}{(u-\alpha)^2} \,du + \frac{\beta-\alpha}{1-\alpha}\int
_{1/\log L}^{\infty} e^{-2\theta c_LL^2\log L   u}e^{-u}\,du,\nonumber
\end{eqnarray}
where the last line uses a change of variable and the results of
Proposition~\ref{prop1locus}. The corresponding quantity when there
are no large events is given by
\[
\beta\int_\beta^1 \frac{e^{-2\theta L^{2u}}}{u^2} \,du + \beta\int
_{1/\log L}^{\infty} e^{-2\theta L^2\log L  u}e^{-u}\,du.
\]
The leading term in each sum is the first one, and we thus see that
if $c_L\ll1$ (i.e., $\rho_L\ll L^{2\alpha}$), then, as expected, the
probability of
IBD is higher in the presence of large events and, moreover, as a
consequence of
shorter genealogies, correlations between gene frequencies
persist over longer spatial scales. See Figure~\ref{fig1locus} for
an illustration
(in which only the leading terms are plotted).
In classical models IBD decays approximately exponentially with the
sampling distance, at least over small scales.
In~\cite{bartonkelleheretheridge2010}, a numerical investigation of
a similar model to that presented here revealed approximately
exponential decay over small scales followed by
a transition to a different exponential rate over somewhat larger scales.
Since the (rigorous) results of Proposition~\ref{prop1locus} only
apply for sufficiently well separated samples, our arguments above
cannot capture this. They do, on
the other hand,
give a clear indication of the reduction of effective population size
due to large events.

\begin{figure}

\includegraphics{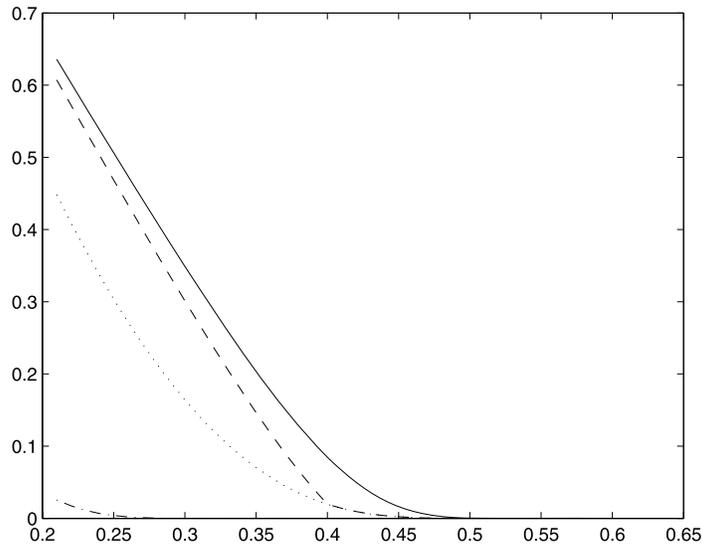}

\caption{Probability of IBD at both loci, as a function of $\beta$.
As in
Figure~\protect\ref{fig1locus}, $L=10^5$, $\alpha=0.1$, $c_L=0.01$ and
$\theta_1=\theta_2=10^{-3}$. The solid line corresponds to the case
$\gamma\geq1$ (complete correlation for any $\beta$), the dotted line
to the case $\gamma\leq\alpha$ (decorrelation for any $\beta$) and the
dashed line to the intermediate case $\gamma=0.4$. The dash--dot line
corresponds to the case without large events, for which $\gamma_L^*$ is
computed from the same parameter values (i.e., $\gamma_L^*=0.2$).}
\label{fig2loci}
\end{figure}

Local bottlenecks are not the only explanations for a reduced effective
population size.
For example, selection or fluctuating population sizes can have the
same effect, and so
we should like to find a more ``personal'' signature of the presence of
demographic events
of different orders of magnitude. The idea that we explore here is to
consider several loci on
the same chromosome, subject to recombination, and to investigate the
pattern of
\textit{linkage disequilibrium} obtained under the assumptions of Section
\ref{results}.
Using the results of Theorems~\ref{theocorrelation} and~\ref{theosemi-correl}, we have
\begin{eqnarray*}
&&\IP_{L^\beta}[\mbox{IBD at both loci}] \\
&&\qquad= \IE_{L^\beta}
\bigl[e^{-2(\theta
_1 \tau_{Aa}^L+\theta_2\tau_{Bb}^L)}\bigr]\\
&&\qquad \approx\IE_{L^\beta}\bigl[e^{-2\theta_1 \tau_{Aa}^L-2\theta
_2\tau
_{Bb}^L}\ind_{\{c_LL^{2\beta}\leq\tau_{Aa}^L=\tau_{Bb}^L\leq
c_LL^{2\gamma}\}}\bigr]\\
&&\qquad\quad{}  + \IE_{L^\beta}\bigl[e^{-2\theta_1 \tau_{Aa}^L }\ind_{\{
\tau
_{Aa}^L>c_LL^{2\gamma}\}}\bigr]\times\IE_{L^\beta}\bigl[e^{-2\theta_2
\tau_{Bb}^L }\ind_{\{\tau_{Bb}^L>c_LL^{2\gamma}\}}\bigr],
\end{eqnarray*}
where $\theta_1$ and $\theta_2$ denote the mutation rates at each locus
and the first integral is $0$ if condition~(\ref{cond}) holds
(i.e., if there is no first period of complete correlation).
By the same computations as in (\ref{eqibd1}), the leading terms in
this expression
are
%
%
\begin{eqnarray}\label{doubleibd}
&&(\beta-\alpha)\int_\beta^\gamma\frac{e^{-2(\theta_1+\theta
_2)c_LL^{2u}}}{(u-\alpha)^2} \,du
\nonumber
\\[-8pt]
\\[-8pt]
\nonumber
&&\qquad{} + (\beta-\alpha)^2\biggl(\int
_\gamma
^1\frac{e^{-2\theta_1c_LL^{2u}}}{(u-\alpha)^2} \,du\biggr)
\biggl(\int
_\gamma^1\frac{e^{-2\theta_2c_LL^{2u}}}{(u-\alpha)^2} \,du\biggr).
\end{eqnarray}
On the other hand, when there are no large events, the analysis of
Lemma~\ref{lemmaseparation}
(with \textit{effective recombination} replaced by \textit{recombination}
and the
separation to attain of the order of $L$)
tells us that the time two lineages initially in the same individual
need to
``decorrelate'' is of the order of $r_L^{-1}\log L$. Here $r_L^{-1}$ is
the expected time
to wait until we see a recombination event, and $\log L$ is (roughly) the
mean number of recombination events before we see one after which the
lineages remain
separated for a duration
$\mathcal{O}(L^t)$ for some $t\in[\beta,1]$.
Hence, when there are only small events, the leading terms in the
probability of
IBD at both loci are
\[
\beta\int_\beta^{\gamma^*_L} \frac{e^{-2(\theta_1+\theta
_2)L^{2u}}}{u^2} \,du + \beta^2\biggl(\int_{\gamma^*_L}^1\frac
{e^{-2\theta
_1L^{2u}}}{u^2} \,du\biggr) \biggl(\int_{\gamma^*_L}^1\frac
{e^{-2\theta
_2L^{2u}}}{u^2} \,du\biggr),
\]
where we have set $\gamma^*_L:=\log(r_L^{-1}\log L)/(2\log L)$ and the
first integral
is again zero if $\beta>\gamma^*_L$. Figure~\ref{fig2loci}
compares the
different curves obtained when (i) we always have decorrelation
($\gamma
\leq\alpha$),
(ii) we always have complete correlation ($\gamma\geq1$), or (iii)
when we have a
transition between these two regimes [$\gamma\in(\alpha,1)$]. As
expected, we see that the
probability of IBD at both loci is higher in the presence of large events
(when $\rho_L\leq L^{2\alpha}$), and there is correlation between the
two loci when
individuals are sampled over large spatial distances.
Furthermore, (\ref{doubleibd}) gives us an idea of how the
correlations between the
two loci decay with sampling distance, as this grows from the radius of
the large
events to the whole population range. Correlations for sampling
distances smaller
than or equal to the size of the large events will be the object of
future work.

%

\section{Genealogies at one locus}\label{section1locus}
In this section we prove Proposition~\ref{prop1locus}. In the
process we introduce a rescaling of the spatial motion of our ancestral
lineages and
collect together several results on the time required to ``gather'' two
lineages to within distance $2R_BL^{\alpha}$ which will also be needed
in Section~\ref{section2loci}.
Since the techniques mirror closely those used in previous work, in the
interests of
brevity, we restrict ourselves
to sketching the proofs and providing references where appropriate.

\textit{Assume for the rest of this section that $\alpha<1$.}

The following local central limit theorem, corresponding to Lemma~5.4
of~\cite{BEV2010},
is the key to understanding the behavior of two lineages. Suppose that
for each $L\in\IN$, $\ell^L$ is a L\'evy process on $\IT(L)$ such that
$\ell^L(1)-\ell^L(0)$ has a covariance matrix of the form $\sigma_L^2
\operatorname{Id}$, and that:
\begin{longlist}[(ii)]
\item[(i)] there exists $\sigma^2>0$ such that $\sigma_L^2\rightarrow
\sigma^2$ as $L\rightarrow\infty$;

\item[(ii)] $\rmE_0[|\ell^L(1)|^4]$ is bounded uniformly in $L$.
\end{longlist}

We shall implicitly suppose that all processes $\ell^L$ are defined on
the same
probability space, and that under the probability measure $\rmP_x$
the L\'evy process we consider starts at $x$. Let $(d_L)_{L\geq1}$ be
a sequence of
positive reals such that $\liminf_{L\rightarrow\infty}d_L>0$ and
$\frac{\log^+(d_L)}{\log L}\rightarrow\eta\in[0,1)$. Finally, let
us write
$p^L(x,t)$ for $\rmP_x[\ell^L(t)\in B(0,d_L)]$ and $\integ{z}$ for
the integer
part of $z\in\IR$.

\renewcommand{\thelemmaa}{\Alph{lemmaa}}
\begin{lemmaa}[(Lemma~5.4 in~\cite{BEV2010})]\label{lemA}
\begin{longlist}
\item[(a)] Let $\e_L:= (\log L)^{-1/2}$. There exists a
constant $C_1<\infty$ such that for every $L\geq2$,
\[
\sup_{t\geq\integ{\e_LL^2}} \sup_{x\in\IT(L)} \frac{\integ
{\e
_LL^2}}{d_L^2} p^L(x,t)\leq C_1.
\]
\item[(b)] If $v_L\rightarrow\infty$ as $L\rightarrow\infty$, then
\[
\lim_{L\rightarrow\infty} \sup_{t\geq\integ{v_LL^2}} \sup
_{x\in\IT
(L)} \frac{L^2}{d_L^2} \biggl| p^L(x,t)-\frac{\pi d_L^2}{L^2}\biggr|=0.
\]
\item[(c)] If $u_L\rightarrow\infty$ as $L\rightarrow\infty$ and
$I(d_L,x):= 1+(|x|^2\vee d_L^2)$, then
\[
\lim_{L\rightarrow\infty} \sup_{x\in\IT(L)} \sup
_{u_LI(d_L,x)\leq
t\leq\e_LL^2} \frac{2\sigma_L^2t}{d_L^2} \biggl|p^L(x,t)-\frac
{d_L^2}{2\sigma_L^2t}\biggr |=0.
\]
\item[(d)] There exists a constant $C_2<\infty$ such that for every
$L\geq1$,
\[
\sup_{t\geq0}\sup_{x\in\IT(L)}\biggl(1+\frac{|x|^2}{d_L^2}
\biggr)p^L(x,t)\leq C_2.
\]
\end{longlist}
\end{lemmaa}

What Lemma~\ref{lemA} shows is that, for times which are large but of order at most
$\cO(L^2)$, $\ell^L$ behaves like two-dimensional Brownian motion (case
$c$), and,
in particular, it has not yet explored the torus enough to ``see'' that
space is limited. On the other hand,
$\ell^L(t)$ is nearly uniformly distributed over $\IT(L)$ at any time
much greater
than $L^2$ (case $b$).

Fix $R>0$. As a direct corollary of this local central limit theorem,
we proved in
Lemma~5.5 of~\cite{BEV2010} that, if $T(R,\ell^L)$ denotes the entrance
time of
$\ell^L$ into the ball $B(0,R)$, then the following inequality holds.

\begin{lemmaa}[(Lemma~5.5 in~\cite{BEV2010})]\label{lemB}
Let $(U_L)_{L\geq1}$ and $(u_L)_{L\geq1}$ be two sequences increasing
to infinity such that $U_LL^{-2}\rightarrow\infty$ as $L\rightarrow
\infty$ and $2u_L\leq L^2(\log L)^{-1/2}$ for every $L\geq1$. Then,
there exist $C_0>0$ and $L_0\in\IN$ such that for every sequence
$(U_L')_{L\geq1}$ satisfying $U_L'\geq U_L$ for each $L$, every $L\geq
L_0$ and all $x\in\IT(L)$,
\[
\rmP_x\bigl[T(R,\ell^L) \in[U_L'-u_L, U_L']\bigr]\leq\frac{C_0 u_L}{L^2}.
\]
\end{lemmaa}

Lemma~\ref{lemB} tells us about the regime in which~$\ell^L$ has
already homogenized
over $\IT(L)$.
Using exactly the same method, but employing parts (c) and (d) of
Lemma~\ref{lemA} rather than
(b), we obtain the analogous result for the regime in which~$\ell^L$
behaves as Brownian
motion on $\IR^2$:
%
\begin{lemma}\label{lemmahitting} If $U_L\leq L^2(\log L)^{-1/2}$ for
each $L\geq1$, $U_L,u_L\rightarrow\infty$ and $u_L/U_L \rightarrow0$
as $L\rightarrow\infty$, then there exist $C_1>0$ and $L_1\in\IN$
such that for every sequence\vadjust{\goodbreak} $(U_L')_{L\geq1}$ satisfying $U_L\leq
U_L'\leq L^2(\log L)^{-1/2}$ for each $L$, for every $L\geq L_1$ and
$x\in\IT(L)$,
\[
\rmP_x\bigl[T(R,\ell^L) \in[U_L'-u_L, U_L']\bigr]\leq\frac{C_1u_L}{U_L'}.\vspace*{-2pt}
\]
\end{lemma}

Let us now introduce the processes to which we wish to apply these results.
For each $L\in\IN$, let $\{\tX_{Aa}^L(t), t\geq0\}$ be the process
recording
the difference between the locations on $\IT(L)$ of the ancestral
lineages of $A$ and $a$
(i.e., the first locus of each of the two individuals sampled).
The process $\tX_{Aa}^L$ is the difference between two dependent
compound Poisson processes.
Under the probability measures we shall use, it is a Markov process
(see Remark~\ref{remmarkovproperty}).
Observe that, because the largest events have radius $R_BL^{\alpha}$,
the lineages have to be within a distance less than
$2R_BL^{\alpha}$ of each other to be hit by the same event. As a
consequence, the law of $\tX_{Aa}^L$ outside $B(0,2R_BL^{\alpha})$ is
equal to that of the difference $\tilde{Y}^L$ of two i.i.d. L\'evy
processes, each of which follows the evolution given in
(\ref{intensities}), and thus is also equal to the law of the motion of
a single lineage run at twice the speed.
We define
the processes $X_{Aa}^L$ and $Y^L$ by
%
%
\begin{equation}\label{scales}
X_{Aa}^L(t)= \frac{1}{L^{\alpha}}  \tX_{Aa}^L(\rho_Lt)
\quad\mbox{and}\quad  Y^L(t)= \frac{1}{L^{\alpha}}  \tilde{Y}^L(\rho_Lt),
\qquad
t\geq0,
\end{equation}
both evolving on $\IT(L^{1-\alpha})$.
Using computations from the proof of Proposition~6.2 in~\cite{BEV2010}
and the
jump intensities given in (\ref{intensities}), we find that
the covariance matrix of $Y^L(1)-Y^L(0)$ is the identity matrix
multiplied by
%
%
\begin{eqnarray}\label{defsigma}\quad
&&2 \biggl\{\frac{u_s\rho_L}{\pi R_s^2L^{2\alpha}} \int_{\IR^2} (x_1)^2
L_{R_s}(x,0) \,dx + \frac{u_B}{\pi R_B^2}\int_{\IR^2} (x_1)^2
L_{R_B}(x,0) \,dx \biggr\}+ o(1)
\nonumber
\\[-8pt]
\\[-8pt]
\nonumber
&&\qquad=: 2\sigma_L^2 + o(1),
\end{eqnarray}
with $\sigma_L^2$ tending to a finite limit $\sigma^2>0$ as
$L\rightarrow\infty$
(by our assumption on $L^{-2\alpha}\rho_L$).
The $o(1)$ remainder here is the error we make by considering $\tilde
{Y}^L$ as evolving on $\IR^2$ instead of $\IT(L)$ (see the proof of
Proposition~6.2 in~\cite{BEV2010}).
Assumption (ii) is also satisfied, and so Lemma~\ref{lemA} and its corollaries
apply to $(Y^L)_{L\geq1}$, with the torus sidelength $L$ replaced by
$L^{1-\alpha}$.
Furthermore, $X_{Aa}^L$ and $Y^L$ follow the same evolution outside
$B(0,2R_B)$ for every $L$.
This will be sufficient to prove Proposition~\ref{prop1locus}: we
shall show that the time the ancestral lineages of $A$ and $a$ need to coalesce
once they are within distance $2R_BL^{\alpha}$ of one another
[or, equivalently, once $X_{Aa}^L$ has entered $B(0,2R_B)$] is negligible
compared to the time they need to be gathered at distance
$2R_BL^{\alpha}$.
It is therefore the ``gathering time'' that dictates the coalescence time
of two lineages starting at separation $|x_L|\gg L^{\alpha}$.

\begin{remark}\label{remmarkovproperty}
It is here that we take advantage of the form of our recombination
mechanism (recall
Remark~\ref{remarksonmodel}).
When $\tX_{Aa}^L(t)\neq0$, its future evolution is determined by the
homogeneous
Poisson point processes of events $\Pi_B^L$ and $\Pi_s^L$, and depends\vadjust{\goodbreak}
only on the
current separation of the two lineages. If $\tX_{Aa}^L(t)=0$, the
situation depends upon
whether the two lineages are in the same individual (i.e., they have
coalesced
and will require a recombination event to separate again), or in two
distinct individuals
at the same spatial location. However, because of the form of our
recombination mechanism, two lineages can jump onto the same location
only if they are descendants of the same parent (in which case they
necessarily coalesce).
This means that provided we choose our initial condition in such a way
that two
lineages in the same spatial location are actually in the same
individual, with
probability one we will never see two lineages in distinct individuals
but the same
spatial location and so
$\tX_{Aa}^L$ is indeed a Markov process under $\IP_{a_L}$.
\end{remark}

\begin{notn}As at the beginning of the section, we assume that all
$Y^L$'s are defined on
the same probability space, and start at $x$ under the probability
measure $\rmP_x$.
Since $X_{Aa}^L$ is a function of the genealogical process of $A,a,B$
and $b$, we retain
the notation $\IP_{a_L}$ when referring to it, and $X_{Aa}^L$ then starts
a.s.~at $L^{-\alpha}x_L$ if $x_L\in\IT(L)$ is the initial separation
between lineages~$A$ and $a$.
\end{notn}

The proof of Proposition~\ref{prop1locus} will require two
subsidiary results.
For each \mbox{$L\in\IN$}, let $T_{Aa}^L$ be the first time the two lineages
$A$ and $a$ are at
separation less than $2R_BL^{\alpha}$. Equivalently, $\rho
_L^{-1}T_{Aa}^L$ is the entrance
time of $X_{Aa}^L$ into $B(0,2R_B)$. By the observation made in the
paragraph preceding
Remark~\ref{remmarkovproperty}, $\rho_L^{-1}T_{Aa}^L$ under $\IP
_{a_L}$ has the same
distribution as $T(2R_B,Y^L)$ under $\rmP_{L^{-\alpha}x_L}$, which
yields the following lemma.
%
\begin{lemma}\label{lemmameeting}
Under the assumptions of Proposition~\ref{prop1locus}, we have
%
%
\begin{eqnarray}
\lim_{L\rightarrow\infty} \IP_{a_L}\bigl[T_{Aa}^L> \rho
_LL^{2(t-\alpha
)}\bigr] &= &\frac{\beta-\alpha}{t-\alpha}\qquad
 \forall t\in[\beta,1]\quad \mbox{and}\label{meeting1} \\
\qquad\lim_{L\rightarrow\infty} \IP_{a_L}\biggl[T_{Aa}^L >\frac
{1-\alpha
}{2\pi\sigma^2} \rho_L L^{2(1-\alpha)}\log L t \biggr] &=& \frac
{\beta
-\alpha}{1-\alpha} e^{-t}\qquad \forall t>0.\label{meeting2}
\end{eqnarray}
Furthermore, for any $\beta_0\in(\alpha,1)$ and $\e>0$, the
convergence in the
first (resp., second) expression is uniform in $\beta,t\in[\beta_0,1]$
(resp., $\beta\in[\beta_0,1]$ and $t\geq\e$) and~$a_L$ such that
$|x_L|\in[L^{\beta}/(\log L),L^{\beta}\log L]$.
\end{lemma}

\begin{pf}When
$\beta=1$, the results
are a weaker version of Proposition~6.2 in~\cite{BEV2010}, in which the
convergence in
(\ref{meeting2}) is uniform over $t\geq0$ and over the set of
sequences $(x_L)_{L\geq1}$
such that $|x_L|\geq L(\log L)^{-1}$ for every $L$. Here, we relax the
condition on
$(x_L)_{L\geq1}$, but since the arguments in the proof of convergence
(without requiring uniformity)
only use the asymptotic behavior of $\log|x_L|$, they are still valid.

If $\beta<1$, the reasoning is the same as in the proofs of Lemma~3.6
in~\cite{ZCD2005}
(note that as above we allow more general sequences of initial
separations at the expense of
the uniformity of the convergence)
and Theorem~2 in~\cite{CD2002}. This does not come as a surprise, since
the same local
central limit theorem applies to both $Y^L$ [on $\IT(L^{1-\alpha})$] and
Z\"ahle, Cox and Durrett's $Y$ [on $\IT(L)\cap\IZ^2$] up to some
constants depending
on the geometry of the geographical patches considered. Hence, since $X_{Aa}^L$
starts from $L^{-\alpha}x_L$ and $\log(L^{-\alpha}|x_L|)/(\log
L)\rightarrow\beta-\alpha$
by assumption, we can write (as in Lemma~3.6 of~\cite{ZCD2005})
\begin{eqnarray*}
&&\lim_{L\rightarrow\infty} \sup_{\beta\leq t \leq\kappa_L}
\biggl|\IP
_{a_L}\bigl[T_{Aa}^L> \rho_LL^{2(t-\alpha)}\bigr]- \frac{\beta
-\alpha
}{t-\alpha}\biggr| \\[-2pt]
& &\qquad= \lim_{L\rightarrow\infty} \sup_{\beta\leq t \leq\kappa
_L}
\biggl|\rmP_{L^{-\alpha}x_L}\bigl[T(2R_B,Y^L)>L^{2(t-\alpha)}\bigr]- \frac
{\beta
-\alpha}{t-\alpha}\biggr| \\[-2pt]
& &\qquad = 0,
\end{eqnarray*}
where $\kappa_L=1-(\log\log L)/(2\log L)$ [so that $L^{2(\kappa
_L-\alpha)}=L^{2(1-\alpha)}/(\log L)$]. Now, as in Lemma~3.8 of \cite
{ZCD2005}, there exists $L_0\in\IN$ and a constant $C$ such that, for
every $L\geq L_0$ and $x\in\IT(L)$,
%
%
\begin{equation}\label{hittingL2}
\rmP_{x}\biggl[Y^L(s)=0 \mbox{ for some }s\in\biggl[\frac
{L^{2(1-\alpha)}}{\log L},L^{2(1-\alpha)}\biggr]\biggr]\leq\frac
{C\log
\log L}{\log L}.
\end{equation}
Combining these two results, we obtain (\ref{meeting1}).

Finally, (\ref{meeting2}) is the analogue of Theorem~2 in~\cite{CD2002}
and can either
be proved using the same technique or in the same way as
Proposition~6.2 in~\cite{BEV2010}
(which, in addition, gives the appropriate constant in the
time-rescaling). The uniform
convergence stated in the second part of Lemma~\ref{lemmameeting} follows
from a direct application of the techniques of~\cite{ZCD2005} and
\cite
{BEV2010} cited above.
\end{pf}

The next result we need is the time that two lineages starting at
separation at most
$2R_BL^{\alpha}$ take to coalesce. Under our assumption that $(\rho
_LL^{-2\alpha})_{L\geq1}$
is bounded, Proposition~6.4(a) in~\cite{BEV2010} applied with $\psi
_L:= L^{\alpha}$
shows that for any sequence $(\phi_L)_{L\geq1}$ tending to infinity,
we have
%
%
\begin{equation}\label{coaltime}
\lim_{L\rightarrow\infty}\sup_{a_L'} \IP_{a_L'}[\tau
_{Aa}^L>\phi
_L\rho_L]=0,
\end{equation}
where the supremum is taken over all configurations $a_L'$ such that
the distance between
the blocks containing $A$ and $a$ is at most $2R_BL^{\alpha}$. Observe
that in
\cite{BEV2010}, only one individual reproduces during an event, and so
if several
lineages are affected by this event, they necessarily coalesce.
Here, the distributions $\lambda_s$ and $\lambda_B$ of the number of
potential parents
are more general, but we assumed that their supports were compact.
Thus, the probability
that several individuals in the area of an event come from the same parent
does not vanish as $L$ tends to infinity, which is all that
we need to prove (\ref{coaltime}).

\begin{remark}
Since (\ref{coaltime}) shows that coming to within $2R_BL^{\alpha}$ is
almost equivalent
to coalescing for two lineages, this is the only point where the distributions
$\lambda_s$ and $\lambda_B$ appear in our discussion.\vadjust{\goodbreak}
\end{remark}

\begin{pf*}{Proof of Proposition~\protect\ref{prop1locus}}
Equipped with these results and the corollaries of Lemma~\ref{lemA}, we can now
write for any given $t\in[\beta,1]$
%
%
\begin{eqnarray}\label{eqcoaltime}
&&\IP_{a_L}\bigl[ \tau_{Aa}^L>\rho_LL^{2(t-\alpha)}\bigr]\nonumber\\
&&\qquad= \IP
_{a_L}\bigl[
\tau_{Aa}^L>\rho_LL^{2(t-\alpha)} ; T_{Aa}^L>\rho_L
\bigl(L^{2(t-\alpha
)}-\log L\bigr)\bigr]
\\
 & &\qquad\quad{} + \IP_{a_L}\bigl[ \tau
_{Aa}^L>\rho
_LL^{2(t-\alpha)} ; T_{Aa}^L\leq\rho_L\bigl(L^{2(t-\alpha)}-\log
L
\bigr)\bigr].\nonumber
\end{eqnarray}
The second term on the right-hand side of (\ref{eqcoaltime}) tends to
zero by the
strong Markov property applied at time $T_{Aa}^L$ and (\ref{coaltime})
with $\phi_L=\log L$.
Then, we have, for each $L$,
\begin{eqnarray*}
&&\bigl|\IP_{a_L}\bigl[ \tau_{Aa}^L>\rho_LL^{2(t-\alpha)};
T_{Aa}^L>\rho
_L\bigl(L^{2(t-\alpha)}-\log L\bigr)\bigr]-\IP_{a_L}
\bigl[T_{Aa}^L>\rho
_LL^{2(t-\alpha)}\bigr]\bigr| \\
&& \qquad\leq\IP_{a_L}\bigl[\rho_L
\bigl(L^{2(t-\alpha)}-\log L\bigr)\leq T_{Aa}^L \leq\rho_LL^{2(t-\alpha
)}
\bigr] \\
&&\qquad= \rmP_{L^{-\alpha}x_L}\bigl[L^{2(t-\alpha)}-\log L \leq T(2R_B,Y^L)
\leq L^{2(t-\alpha)}\bigr],
\end{eqnarray*}
which tends to zero by Lemma~\ref{lemmahitting} applied with $L$
replaced by
$L^{1-\alpha}$ (the size of the torus on which $Y^L$ evolves) if $t<1$,
and by
(\ref{hittingL2}) if $t=1$. Lemma~\ref{lemmameeting} enables us to
deduce (a).

For (b), the same technique applies but with the last argument
replaced by the use of Lemma~\ref{lemB}.
\end{pf*}

Proposition~\ref{prop1locus} is, in fact, a particular case of a
more general result
which we shall use in Section~\ref{subsproofs} (with $k=4$). Suppose
we follow the ancestry
at one locus of $k\geq2$ different individuals. By analogy with above,
we label individuals $1,\ldots,k$, we write $x_{ij}^L$ for the initial
separation of
lineages $i$ and $j$, $T_{ij}^L$ for the time at which their ancestral
lineages first come
within $2R_BL^{\alpha}$ and $\tau_{ij}^L$ for their coalescence time.
We also write
$T_*^L$ (resp., $\tau_*^L$) for the minimum over $\{i\neq j\}$ of the
$T_{ij}^L$'s
(resp., the $\tau_{ij}^L$'s).
Although (in the same way as above) we could state a result for a more
general sequence
$(a_L)_{L\geq1}$ of initial configurations,
for the proof of Theorem~\ref{theocorrelation} we shall need some
uniformity in
the convergence. For this reason, we consider $\Gamma(L,k,\eta)$, the
set of all
configurations of $k$ lineages on $\IT(L)$ such that all pairwise
distances $|x_{ij}^L|$
belong to $[L^{\eta}/(\log L),L^{\eta}\log L]$.
%
\begin{propn}\label{propkingman}For any $\beta\in(\alpha,1]$, $\e>0$
and $i\neq j$, we have \setlength\arraycolsep{0pt}
\begin{eqnarray*}
&&\lim_{L\rightarrow\infty}\sup_{\beta\leq\eta\leq t\leq1}\sup
_{a_L\in\Gamma(L,k,\eta)}\biggl|\IP_{a_L}\bigl[\tau_*^L=\tau
_{ij}^L\leq
\rho_LL^{2(t-\alpha)}\bigr]-\frac{1}{{k\choose 2}}\biggl(1-
\biggl(\frac{\eta
-\alpha}{t-\alpha}\biggr)^{{k}\choose{2}}\biggr)\biggr|\\
&&\qquad=0,\\
&&\lim_{L\rightarrow\infty}\sup_{t\geq\e,\beta\leq\eta\leq
1}\sup
_{a_L\in\Gamma(L,k,\eta)}\biggl|\IP_{a_L}\biggl[\tau_*^L=\tau
_{ij}^L\leq
\frac{1-\alpha}{2\pi\sigma^2} \rho_LL^{2(1-\alpha)}\log L t
\biggr]\\
& &\hspace*{150pt}\qquad{}
-\frac{1}{{k \choose 2}}\biggl(1-\biggl(\frac{\eta-\alpha}{1-\alpha
} e^{-t}\biggr)^{{k}\choose{2}}\biggr)\biggr|= 0.
\end{eqnarray*}
The same is true with $\tau^L$ replaced by $T^L$.
\end{propn}

In essence, Proposition~\ref{propkingman} tells us that on the timescale
$\rho_LL^{2(t-\alpha)}, t\in[\eta,1]$,
the time of the first coalescence (or
of the first ``gathering'') is approximately the same as that of the
first merger in a Kingman
coalescent timechanged by
$\log(\frac{t-\alpha}{\eta-\alpha})$, and that the
approximation is uniform
over $\eta$'s bounded away from $\alpha$. Moreover, asymptotically,
just as in the
Kingman coalescent, each pair of lineages has the same chance to be the
first to coalesce.
On the other hand, on the timescale $\frac{1-\alpha}{2\pi\sigma
^2}\rho
_LL^{2(1-\alpha)}\log L  t$,
conditional on $T_*^L>\rho_LL^{2(1-\alpha)}$, the asymptotic behavior
corresponds to Kingman's coalescent run at speed $1$.

\begin{pf*}{Sketch of proof}
The proof of Proposition~\ref{propkingman} is a straightforward
adaptation of those of
Lemma~4.2 and of Lemma~5.2 in~\cite{ZCD2005} (see also the comments
given in the paragraph
following the proof of Lemma~4.2). The interested reader will also find there
references to earlier results for the random walks with instantaneous
coalescence
which are dual to the two-dimensional voter model.
\end{pf*}

Let us end this section by recalling a lemma of~\cite{BEV2010} and by
stating an
analogous result. For every $L\in\IN$, $i\neq j$ and $t\geq0$, let
$\tilde{X}_{ij}^L(t)$ be
the separation [on $\IT(L)$ at time $t$] of lineages $i$ and $j$.

\begin{lemmaa}[(Lemma~6.9 in~\cite{BEV2010})]\label{lemC}
Suppose
$k=4$ and
%
%
\begin{equation}\label{initialsep}
\lim_{L\rightarrow\infty}\frac{\min_{i\neq j}\log|x_{ij}^L|}{\log
L}=\lim_{L\rightarrow\infty}\frac{\max_{i\neq j}\log
|x_{ij}^L|}{\log L}=1.
\end{equation}
Then,
\begin{eqnarray*}
\lim_{L\rightarrow\infty}\IP_{a_L}\biggl[\tau_*^L=\tau_{12}^L ;
|\tilde{X}_{13}^L(\tau_*^L)|\leq\frac{L}{\log L}\biggr]&=&0, \\
\lim_{L\rightarrow\infty}\IP_{a_L}\biggl[\tau_*^L=\tau_{12}^L ;
|\tilde{X}_{34}^L(\tau_*^L)|\leq\frac{L}{\log L}\biggr]&=&0.
\end{eqnarray*}
These results are also true if $\tau^L$ is replaced by $T^L$.
\end{lemmaa}

In words, when two lineages meet and coalesce, with probability tending
to one the others
are at distance at least $L/\log L$ of each other and of the coalescing pair
(in particular, such a merger involves at most two lineages at a time).
When the
initial distance between the lineages is of the order of $L^{\beta}$
with $\beta<1$,
we have instead:
%
\begin{lemma}\label{lemmasep}
Suppose again $k=4$ and the limit in (\ref{initialsep}) is equal to
$\beta\in(\alpha,1)$. Then,
\begin{eqnarray*}
&&\lim_{L\rightarrow\infty}\IP_{a_L}\biggl[\tau_*^L=\tau_{12}^L\leq
\frac
{\rho_LL^{2(1-\alpha)}}{\log L} ; |\tilde{X}_{13}^L(\tau
_*^L)
|\notin\biggl[\frac{L^{\alpha}}{\log L} \frac{\sqrt{\tau
_*^L}}{\sqrt
{\rho_L}},L^{\alpha}\log L \frac{\sqrt{\tau_*^L}}{\sqrt{\rho
_L}}
\biggr]\biggr]\\
&&\qquad =0, \\
&&\lim_{L\rightarrow\infty}\IP_{a_L}\biggl[\tau_*^L=\tau_{12}^L\leq
\frac
{\rho_LL^{2(1-\alpha)}}{\log L} ; |\tilde{X}_{34}^L(\tau
_*^L)
|\notin\biggl[\frac{L^{\alpha}}{\log L} \frac{\sqrt{\tau
_*^L}}{\sqrt
{\rho_L}},L^{\alpha}\log L \frac{\sqrt{\tau_*^L}}{\sqrt{\rho
_L}}
\biggr]\biggr]\\
&&\qquad =0.
\end{eqnarray*}
The result is also true if $\tau^L$ is replaced by $T^L$.
\end{lemma}

Notice the rescalings of time by $\rho_L$ and space by $L^{\alpha}$
introduced in (\ref{scales})
under which the behavior of the lineages is close to that of finite
variance random walks.
In fact, although their formulations are rather different, Lemma~\ref
{lemmasep} is very
similar to Lemma~1 in~\cite{CG1986} or Lemma~5.1 in~\cite{ZCD2005} for
coalescing random walks.

\begin{pf*}{Sketch of proof of Lemma~\protect\ref{lemmasep}}
The method of proof is identical to that of Lemma~6.9 in \cite
{BEV2010}, to which we refer for
more complete arguments. It is based on two facts. First, by time
$\rho_LL^{2(1-\alpha)}/(\log L)$ the separation of the lineages is
never on the order of
the side of the torus. Second, if $T_*^L=T_{12}^L$, then
$L^{-\alpha}\tilde{X}_{13}^L(\rho_L\cdot)$ and $L^{-\alpha}\tilde
{X}_{34}^L(\rho_L\cdot)$,
considered separately, follow the same law as the difference of two
independent lineages
(on $\IR^2$, by the first fact) conditioned on not entering $B(0,2R_B)$
before $T_*^L/\rho_L$.
By Lemma~\ref{lemmameeting}, with high probability $T_*^L/\rho_L
\gg
L^{2(\beta-\alpha)}$,
and so the result for $T^L$ follows from a standard central limit theorem.

The modifications needed for $\tau^L$ use the very rapid coalescence of
two lineages
gathered at distance $2R_BL^{\alpha}$ to obtain that, with probability
tending to $1$,
if $\tau_*^L=\tau_{12}^L$, then no other pairs of lineages come within
$2R_BL^{\alpha}$ of one another before time~$\tau_*^L$.
An application of Lemma~\ref{lemmasep} (with $T^L$) completes the
proof.
\end{pf*}

\section{Genealogies at two loci}\label{section2loci}
From now on, we work with the rescaling of time and space introduced in
(\ref{scales}).
As we saw in the previous section, these are the appropriate scales on which
to understand the behavior of a collection of independent processes following
the dynamics driven by (\ref{intensities}). Because our lineages move
independently as
long as they are at distance greater than $2R_B$ (in rescaled units) of
each other,
it is also the relevant regime in which to understand ``gathering'' and
coalescence of ancestral lineages.

The aim of Sections~\ref{subseffectiverecomb} and~\ref{subsdecorrelation} is to understand
how two lineages, initially present in the same individual, can
``decorrelate'' and how much
time they need to do so. Once this phenomenon is understood for two
lineages, we can
consider the more complex situation described in the Section \ref
{results} and prove
Theorems~\ref{theocorrelation} and~\ref{theosemi-correl}. This is
achieved in Section~\ref{subsproofs}.

\subsection{Effective recombination time}\label{subseffectiverecomb}
For every $L$, let $X_{AB}^L$ be the process that records the
(rescaled) difference
between the locations of the lineages labeled $A$ and $B$. Recall that
under our
working assumptions, these lineages start within the same individual
(in other
words, $A$ and $B$ belong to the same block of the marked partition $a_L$).

By construction, recombination occurs only during small events.
In our rescaled space and time units, a recombination event results in
a separation of the
lineages of $\cO(L^{-\alpha})$, and then
small events affect them at rate $\cO(\rho_L)$.
Hence, it is very likely that (in our rescaled time units) the lineages
very rapidly coalesce and have to wait for the next recombination event
[i.e., roughly $(\rho_L r_L)^{-1}$ units of rescaled time] to be
geographically
separated again, and so on. An efficient way for the lineages to escape this
``flickering'' due to small events is for a large event to send them to a
separation
of $\cO(1)$. This necessarily occurs at a time when $X_{AB}^L\neq0$.
Thus, let us define $S^L$ as the first time $t$ at which at least one
of the two
lineages is affected by a large event and $X_{AB}^L(t-)\neq0$ [which
does not prohibit
$X_{AB}^L(t)=0$]. We call $S^L$ the \textit{effective recombination time}.
Its large-$L$ behavior is given by the following proposition.

\begin{propn}\label{propeffectiverecomb}
There exist $\theta_1, \theta_2>0$ such that for every $\theta
>\theta
_2$ and every nonvanishing sequence $(\phi_L)_{L\geq1}$ satisfying
$\phi_L\leq L^2/(\rho_L\log L)$ for every $L$, we have for $L$ large enough
\[
\IP_{a_L}\biggl[S^L\geq\phi_L\biggl(1+\frac{\theta\log(\phi
_L\rho
_L)}{r_L\rho_L}\biggr)\biggr]\leq e^{-\theta_1\phi_L}+e^{-(\theta
-\theta
_2)\phi_L\log(\phi_L\rho_L)}.
\]
\end{propn}

The idea of the proof of Proposition~\ref{propeffectiverecomb} is to
show that, with very high probability, the number of visits to $0$ of
$X_{AB}^L$ before
it has accumulated a time $\phi_L$ outside $0$
is less than $\phi_L\log(\phi_L\rho_L)$. Since each visit lasts a time
proportional
to $(r_L\rho_L)^{-1}$, the total amount of time it takes for $X_{AB}^L$
to accumulate
$\phi_L$ units of time outside zero is at most of the order of
$\phi_L+\phi_L\log(\phi_L\rho_L)/(r_L\rho_L)$. The probability
that by
this time the
two lineages have not been affected by a large event while in distinct locations
is bounded by a quantity of the form $e^{-C\phi_L}$.

Let us write $\mathcal{R}_L(x)$ for the rate at which at least one of
the lineages is
affected by a large event when $X_{AB}^L=x$, and recall that time is
rescaled by a factor $\rho_L$. From the expression for the intensity of
$\Pi_B^L$, we can find a constant $C_B>0$ such that $\mathcal
{R}_L(x)\geq C_B$ for
all $x\in\IT(L^{1-\alpha})\setminus\{0\}$ [in fact, one can even show
that the function $x\mapsto\mathcal{R}_L(x)$ is increasing in $|x|$,
and so one can take $C_B:= \mathcal{R}_L(0)>0$]. Let $\hX^L$ be a
$\IT
(L^{1-\alpha})$-valued Markov process distributed in the same way as
the difference between two lineages subject only to the events of $\Pi
_s^L$, and $\hS^L$
be an exponential random variable with instantaneous rate
$\mathcal{R}_L(\hX^L(t))\ind_{\{\hX^L(t)\neq0\}}$. By the
preceding remark,
$\hS^L$ is stochastically bounded by an exponential random variable
with instantaneous
rate $C_B\ind_{\{\hX^L(t)\neq0\}}$. Because large events have no
effect when $X_{AB}^L=0$,
the law of the stopped process $\{X_{AB}^L(t), t\in[0,S^L]\}$ is the
same as that
of $\{\hX^L(t), t\in[0,\hS^L]\}$. Thus, for the proof of
Proposition~\ref{propeffectiverecomb}
we work with $\hX^L$ and $\hS^L$ and use $\rmP_x$ to denote the law of
$\hX^L$ under which
$\IP[\hX^L(0)=x]=1$.

For each $L\in\IN$, let us define the stopping times $(\hQ
_i^L)_{i\geq
0}$ and $(\hat{q}_i^L)_{i\geq0}$ by $\hQ_0^L=\hat{q}_0^L=0$ and for every
$i\geq1$,
\begin{eqnarray*}
\hQ_i^L& := & \inf\{t\geq\hat{q}_{i-1}^L\dvtx  \hX^L(t)\neq0
\} ,\\
\hat{q}_i^L &:= & \inf\{t\geq\hQ_i^L\dvtx  \hX^L(t)=0\}.
\end{eqnarray*}
[Note that $\hat{Q}_1^L=0$ if $\hat{X}^L(0)\neq0$, in which case
$\hat
{q}_1^L$ is the first hitting time of $0$.]
By construction, the random variables $(\hQ_i^L-\hat
{q}_{i-1}^L)_{i\in\IN
}$ are i.i.d.
and distributed according to an exponential random variable with parameter
$C_{\mathrm{rec}}r_L\rho_L$, where $C_{\mathrm{rec}}:= \pi R_s^2u_s
(1-\lambda_s(\{1\}))>0$
(the last factor arises since the number of reproducing individuals
needs to be
greater than one for recombination to occur).
We have the following result
for the excursions of $\hX^L$ away from $0$.
%
\begin{lemma}\label{lemmaoutside0}
There exist $C_e>0$ and $u_e>0$ such that for every $L\geq1$ and
$u_e\leq u\leq L^2/(\log L)$, for every $x\in B(0,2R_sL^{-\alpha
})\setminus\{0\}$,
\[
\rmP_x[\hat{q}_1^L>u\rho_L^{-1} ]\geq\frac{C_e}{\log u}.
\]
\end{lemma}

\begin{pf}
Here (and only
here) it is
easier to work with the initial time and space units and show that the
probability
of an excursion outside $0$ of length greater than~$u$ is bounded from
below by $C_e/(\log u)$ when $u$ is large. Let us thus define $\tX^L$ by
$\tX^L(t):= L^{\alpha}\hX^L(\rho_L^{-1}t)$ for all $t\geq0$, with the
understanding that $\tX^L$ starts at $L^{\alpha}x$ under the probability
measure $\rmP_{x}$.

The desired result is shown in~\cite{RR1966} for standard discrete
space random walks
whose jumps have finite variance as well as for Brownian motion (with
the hitting
time of $0$ replaced by the entrance time into a ball of fixed radius)
in two dimensions.
To see why it is true for $\tX^L$ on $\IT(L)$, observe first that by
time $L^2/(\log L)$,
the process $\tX^L$ does not see that space is limited, and so it
behaves as though
it were moving in $\IR^2$.
More precisely, there exists a constant $C>0$ such that for all $z\in
B(0,6R_sL^{-\alpha})$,
\[
\rmP_z\biggl[\sup_{u\leq L^2/(\log L)}|\tX^L(u)|> \frac
{L}{3}\biggr]\leq\frac{C}{\log L}.
\]
(Use the $L^2$-maximal inequality and the fact that $|\tX^L|$ is
bounded by the
corresponding quantity for the same process defined on $\IR^2$, which is
proportional to $L^2/(\log L)$ by equation (22) in~\cite{BEV2010}).
Hence, let us
assume that $\tX^L$ is defined on $\IR^2$ instead of $\IT(L)$. Since
the evolution due to
small events depends on $L$ only through the torus sidelength, with our
new convention
all $\tX^L$'s have the same distribution and we can drop the exponent
$L$ in the notation.
For the same reason, we also write $\tilde{q}_1$ for the random times
$\rho_L\hat{q}_1^L$,
that is, the length of the first excursion outside $0$ of $\tX$.

Let $\tilde{\rmT}_{(4R_s)}$ denote the first time $\tX$ leaves
$B(0,4R_s)$ [and
so $\tX(\tilde{\rmT}_{(4R_s)})\in B(0,6R_s)\setminus B(0,4R_s)$ by our
assumption on the
jump sizes], and let $\tilde{\rmT}_{[2R_s]}$ be the first return time
of $\tX$
into $B(0,2R_s)$ after $\tilde{\rmT}_{(4R_s)}$. We have for every
$x\in
B(0,2R_sL^{-\alpha})\setminus\{0\}$,
%
%
\begin{eqnarray}\label{twoinfs}
&&\rmP_x[\tilde{q}_1>u]\nonumber\hspace*{-15pt}\\
&&\qquad\geq \rmP_x\bigl[\tilde{q}_1>u ; \tilde
{\rmT}_{(4R_s)}<\tilde{q}_1\bigr]\nonumber\hspace*{-15pt}\\
&&\qquad\geq \rmP_x\bigl[\tilde{q}_1-\tilde{\rmT}_{(4R_s)}>u; \tilde
{\rmT
}_{(4R_s)}< \tilde{q}_1 \bigr]
\nonumber\hspace*{-15pt}
\\[-8pt]
\\[-8pt]
\nonumber
&& \qquad= \rmE_x\Bigl[\ind_{\{\tilde{\rmT}_{(4R_s)}<\tilde{q}_1\}}\rmP
_{\tX
(\tilde{\rmT}_{(4R_s)})}[\tilde{q}_1>u]\Bigr]\hspace*{-15pt} \\
&&\qquad\geq \rmE_x\bigl[\ind_{\{\tilde{\rmT}_{(4R_s)}<\tilde{q}_1\}
}\rmP_{\tX
(\tilde{\rmT}_{(4R_s)})}\bigl[\tilde{\rmT}_{[2R_s]}>u\bigr]\bigr]
\nonumber\hspace*{-15pt}
\\
&&\qquad\geq \Bigl(\inf_{B(0,2R_s)\setminus\{0\}}\rmP_{L^{-\alpha
}y}
\bigl[\tilde{\rmT}_{(4R_s)}<\tilde{q}_1\bigr]\Bigr)\Bigl(\inf
_{B(0,6R_s)\setminus B(0,4R_s)}\rmP_{L^{-\alpha}z}\bigl[\tilde{\rmT
}_{[2R_s]}>u\bigr]\Bigr).\nonumber\hspace*{-15pt}
\end{eqnarray}
The first infimum is strictly positive. To see this, note that
$\rmP_{L^{-\alpha}y}[\tilde{\rmT}_{(4R_s)}<\tilde{q}_1]$ is bounded
from below by the
probability that the first four small events affecting the lineages
send them to a
distance at least $4R_s$ of each other before they coalesce, and the
infimum over
$B(0,2R_s)\setminus\{0\}$ of the latter probability is positive since
$u_s<1$ (if $u_s=1$,
only one of the lineages can be in the geographical range of such
separating events, and
so their probability of occurrence shrinks to $0$ as $|y|\rightarrow0$).

For the second infimum in (\ref{twoinfs}), we use the same
construction as in the proof
of Skorokhod embedding (see, e.g.,~\cite{BIL1995}) to write the path of
$\tX$ as that of a
standard Brownian motion $W$ considered at particular times. More precisely,
if $(\tilde{\sigma}_i)_{i\in\IN}$ is the sequence of jump times of
$\tX
$, we can
find a sequence of Brownian stopping times $(\sigma_i)_{i\in\IN}$
such that
$(W(\sigma_i))_{i\geq0}$ has the same joint distributions as
$(\tX(\tilde{\sigma}_i))_{i\geq0}$. For every $i\in\IN$,
conditional on
$W(\sigma_{i-1})$, $\sigma_i$ is the first time greater than $\sigma
_{i-1}$ at which $W$
leaves $B(W(\sigma_{i-1}),l_i)$, where the random variable $l_i$ is
independent of $W$ and
of $\{\sigma_j,j<i\}$ and has the same distribution as the length of
the first jump of
$\tX$. As a consequence, if $\tilde{n}(u):= \max\{i\dvtx  \tilde{\sigma
}_i\leq u\}$,
by comparing the paths of $\tX$ and of $W$ we obtain
\[
\rmP_{L^{-\alpha}z}\bigl[\tilde{\rmT}_{[2R_s]}>u\bigr]\geq\rmP
_{z}
\bigl[W(t)\notin B(0,2R_s), \forall t\leq\sigma_{\tilde{n}(u)}\bigr].
\]
Now, each $\tilde{\sigma}_i-\tilde{\sigma}_{i-1}$ is stochastically
bounded from below
by an exponential random variable with positive parameter $k_1>0$, and
so by standard large
deviation results we can find $k_2>0$ large enough and $k_3>0$ such
that for all $u>1$ and
$y\in\IR^2$,
\[
\rmP_y[\tilde{n}(u)>k_2 u]\leq e^{-k_3 u}.
\]
By construction, each $\sigma_i-\sigma_{i-1}$ is stochastically bounded
from above by
the first time Brownian motion started at $0$ leaves $B(0,2R_s)$, which
also has an
exponential moment. Hence, there exist $k_4,k_5>0$ such that for all
$u>1$ and $y\in\IR^2$,
\[
\rmP_y\bigl[\sigma_{\lfloor k_2 u\rfloor+1}>k_4 u\bigr]\leq e^{-k_5 u}.
\]
Using these bounds and the result already established in~\cite{RR1966}
for Brownian motion
at time $k_4u$, Lemma~\ref{lemmaoutside0} is proved.
\end{pf}

We now have all the ingredients we require to prove Proposition~\ref
{propeffectiverecomb}.

\begin{pf*}{Proof of Proposition~\protect\ref{propeffectiverecomb}}
Set
%
%
\begin{equation}\label{defpsiL}
\psi_L:= \phi_L\biggl(1+\frac{\theta\log(\phi_L\rho_L)}{r_L\rho
_L}\biggr),
\end{equation}
and call $\mathbf{t}(\psi_L)$ the time $\hX^L$ spends away from $0$
before time $\psi_L$.
We have, for every $L$,
\begin{eqnarray*}
\rmP_0[\hS^L\geq\psi_L]&=& \rmP_0[\hS^L\geq\psi
_L; \mathbf{t}(\psi_L)\leq\phi_L ] + \rmP_0[\hS^L\geq\psi
_L; \mathbf{t}(\psi_L)> \phi_L]\\
&\leq& \rmP_0[\mathbf{t}(\psi_L)\leq\phi_L]
+e^{-C_B\phi_L},
\end{eqnarray*}
where $C_B$ is the lower bound on the rate of effective large events
introduced just below the
statement of the proposition. Next, if we set $\hat{k}_L:= \sup \{
i\dvtx  \hQ_i^L\leq\psi_L\}$,
that is, $\hat{k}_L$ is the number of excursions of $\hX^L$ away from
$0$ which start before time~$\psi_L$, we can write
\begin{eqnarray*}
\rmP_0[\mathbf{t}(\psi_L)\leq\phi_L]&=& \rmP_0
[\mathbf
{t}(\psi_L)\leq\phi_L; \hat{k}_L\leq\phi_L \log(\phi_L\rho
_L)] \\
& &{} + \rmP_0[\mathbf{t}(\psi_L)\leq\phi_L; \hat{k}_L>\phi
_L \log
(\phi_L\rho_L)].
\end{eqnarray*}
On the one hand,
\begin{eqnarray*}
\rmP_0[\mathbf{t}(\psi_L)\leq\phi_L; \hat{k}_L>\phi_L \log
(\phi
_L\rho_L)]&\leq& \rmP_0\Biggl[\sum_{i=1}^{\lfloor\phi_L\log
(\phi
_L\rho_L)\rfloor}(\hat{q}_i^L-\hQ_i^L)\leq\phi_L\Biggr] \\
&\leq& \rmP_0[\hat{q}_1^L-\hQ_1^L\leq\phi_L]^{\lfloor
\phi_L \log
(\phi_L\rho_L)\rfloor} \\
&\leq& \biggl(1-\frac{C_e}{\log(\phi_L\rho_L)}\biggr)^{\lfloor
\phi_L
\log(\phi_L\rho_L)\rfloor} \\
&\leq& e^{-C_e' \phi_L},
\end{eqnarray*}
for a constant $C_e'>0$ and $L$ large enough. The second line is
obtained by an obvious recursion using the strong Markov property at
the successive times $\hat{q}_i^L$ in decreasing order, and the third line
uses Lemma~\ref{lemmaoutside0} [recall that by assumption on $\phi
_L$, we have $\phi_L\rho_L\rightarrow\infty$ and $\phi_L\rho
_L\leq
L^2/(\log L)$]. Hence, we can set $\theta_1:= C_B\wedge C_e'$. On the
other hand,
\begin{eqnarray*}
&&\rmP_0[\mathbf{t}(\psi_L)\leq\phi_L; \hat{k}_L\leq\phi_L
\log
(\phi_L\rho_L)]\\
& &\qquad\leq\rmP_0\Biggl[\sum_{i=1}^{\lfloor\phi_L\log(\phi_L\rho
_L)\rfloor+1}(\hQ_i^L-\hat{q}_{i-1}^L)\geq\psi_L-\phi_L\Biggr] \\
& &\qquad= \rmP_0\Biggl[\exp\Biggl\{r_L\rho_L \sum_{i=1}^{\lfloor\phi
_L\log
(\phi_L\rho_L)\rfloor+1}(\hQ_i^L-\hat{q}_{i-1}^L)\Biggr\}\geq\exp
\{
r_L\rho_L(\psi_L-\phi_L)\}\Biggr] \\
& &\qquad\leq e^{-\theta\phi_L\log(\phi_L \rho_L)}\rmE_0\Biggl[\exp
\Biggl\{
r_L\rho_L \sum_{i=1}^{\lfloor\phi_L\log(\phi_L\rho_L)\rfloor
+1}(\hQ
_i^L-\hat{q}_{i-1}^L)\Biggr\}\Biggr],
\end{eqnarray*}
where the last line uses the Markov inequality. As we pointed out
above, the random
variables $r_L\rho_L(\hQ_i^L-\hat{q}_{i-1}^L)$ are i.i.d. with law
$\operatorname
{Exp}(C_{\mathrm{rec}})$.
Therefore, we can write for a constant $\theta_2>0$
\[
\rmP_0[\mathbf{t}(\psi_L)\leq\phi_L; \hat{k}_L\leq\phi_L
\log
(\phi_L\rho_L)]\leq e^{-(\theta-\theta_2)\phi_L\log(\phi_L
\rho_L)}.
\]
Combining these results, the proof of Proposition~\ref{propeffectiverecomb} is complete.
\end{pf*}

Finally, let us use Proposition~\ref{propeffectiverecomb} to obtain
some estimates on
the time two lineages starting in the same individual need to reach a
separation at which
they start to evolve independently. The following lemma will be a key
result for the proof
of Proposition~\ref{propdecorr} in the next section. For every $L\in
\IN$,
let $\rmT_{(3R_B)}^L$ denote the exit time of $X_{AB}^L$ from $B(0,3R_B)$.
%
\begin{lemma}\label{lemmaseparation}
There exists a constant $\theta_3>0$ such that if $(\phi_L)_{L\geq1}$
is as in Proposition~\ref{propeffectiverecomb}, $\phi_L\rightarrow
\infty$ as $L\rightarrow\infty$ and $\theta>\theta_2$, there exists
$L_0=L_0(\theta, (\phi_L)_{L\in\IN})$ such that for every $L\geq L_0$,
\[
\IP_{a_L}\biggl[\rmT_{(3R_B)}^L\geq\phi_L\biggl(1+\frac{\theta
\log(\phi
_L\rho_L)}{r_L\rho_L}\biggr)\biggr]\leq\sqrt{\phi_L}  e^{-\theta
_3 \sqrt
{\phi_L}}.
\]
\end{lemma}

\begin{pf}
For conciseness, we again use the notation $\psi_L$ introduced in~(\ref{defpsiL}).
This time we define $Q_0^L=q_0^L=0$ and
\begin{eqnarray*}
Q_i^L& := &\inf\{t> q_{i-1}^L\dvtx  t \mbox{ is the epoch of an effective recombination}\}, \\
q_i^L& := &\inf\{t\geq Q_i^L\dvtx  X_{AB}^L(t)=0 \mbox{ or } X_{AB}^L(t)\notin B(0,3R_B)\},\\
k_L& := & \max\bigl\{i\dvtx   Q_i^L\leq\rmT_{(3R_B)}^L\bigr\}.
\end{eqnarray*}

First, we claim that there exists a constant $\wp>0$ independent of $L$
such that, for~$L$ large enough,
$k_L+1$ is stochastically bounded by a geometric random variable with
success probability
$\wp$. In other words, the probability that $X_{AB}^L$ starting at
$x\in B(0,3R_B)\setminus\{0\}$ leaves $B(0,3R_B)$ before hitting $0$
is bounded from
below by $\wp$, independently of $x$. The proof of this claim is given
in the first
paragraph of the proof of Lemma~6.6 in~\cite{BEV2010}. (The quantity
$\wp$ is taken to be
the probability that a sequence of large events sends the lineages to a
distance of at
least $3R_B$ without meanwhile being counteracted by small events
bringing them too close
together.) As a consequence, for any large $L$,
\[
\IP_{a_L} \bigl[\rmT_{(3R_B)}^L\geq\psi_L\bigr] \leq\IP_{a_L}
\bigl[\rmT
_{(3R_B)}^L\geq\psi_L  ; k_L< \sqrt{\phi_L}\bigr]+ (1-\wp
)^{\sqrt
{\phi_L}}.
\]
Next, let us write \setlength\arraycolsep{1pt}
%
%
\begin{eqnarray} \label{effect1}
&&\IP_{a_L} \bigl[ \rmT_{(3R_B)}^L\geq\psi_L ; k_L<\sqrt{\phi
_L}
\bigr]
\nonumber
\\[-8pt]
\\[-8pt]
\nonumber
& &\qquad= \IP_{a_L}\Biggl[\rmT_{(3R_B)}^L\geq\psi_L ; k_L<\sqrt{\phi
_L} ; \sum_{i=1}^{k_L}(Q_i^L-q_{i-1}^L)\geq\frac{\psi_L}{2}\Biggr]\\
& &\qquad\quad{}   + \IP_{a_L}\Biggl[\rmT_{(3R_B)}^L\geq\psi_L ; k_L<\sqrt
{\phi
_L} ; \sum_{i=1}^{k_L}(Q_i^L-q_{i-1}^L)< \frac{\psi_L}{2}
\Biggr].\label{effect2}
\end{eqnarray}
The quantity in (\ref{effect1}) is bounded by
%
%
\begin{eqnarray}
\label{boundexcursion}
\IP_{a_L}\Biggl[\sum_{i=1}^{\lfloor\sqrt{\phi_L}\rfloor}
(Q_i^L-q_{i-1}^L)\geq\frac{\psi_L}{2}\Biggr] & = & 1- \IP
_{a_L}
\Biggl[\sum_{i=1}^{\lfloor\sqrt{\phi_L}\rfloor}(Q_i^L-q_{i-1}^L)<
\frac{\psi_L}{2}\Biggr]\nonumber\\
&\leq& 1-\IP_{a_L}\biggl[\forall i\leq\lfloor\sqrt{\phi
_L}\rfloor, Q_i^L-q_{i-1}^L< \frac{\psi_L}{2 \sqrt{\phi_L}}\biggr]\\
& \leq& 1-\biggl(1-\sup_{a_L'} \IP_{a_L'}\biggl[S^L \geq\frac
{\psi
_L}{2 \sqrt{\phi_L}}\biggr]\biggr)^{\lfloor\sqrt{\phi_L}\rfloor},\nonumber
\end{eqnarray}
where the last line is obtained by recursion (notice that,
conditionally on $q_{i-1}^L$, $Q_i^L-q_{i-1}^L$ has the same law as the
effective recombination time $S_L$) and the supremum is taken over all initial
configurations $a_L'$ in which lineages $A$ and $B$ are either at
distance $0$ or at distance
greater than $3R_B$. We can in fact restrict our attention to the set
of configurations in
which $A$ and $B$ belong to the same block. Indeed, if
$|X_{AB}^L(0)|>3R_B$, we can decompose
the probability that $S^L\geq\psi_L/(2 \sqrt{\phi_L})$ into the sum of:
\begin{itemize}
\item the probability that $S^L\geq\psi_L/(2\sqrt{\phi_L})$ and
$X_{AB}^L$ does not hit $0$ before time $\psi_L/(4\sqrt{\phi_L})$,
which decreases like $e^{-C\psi_L/\sqrt{\phi_L}}$ since the rate at
which large events affect the lineages when $X_{AB}^L\neq0$ is bounded
from below by a positive constant;
\item the probability that $S^L\geq\psi_L/(2 \sqrt{\phi_L})$ and
$X_{AB}^L$ hits $0$ before time $\psi_L/\break(4\sqrt{\phi_L})$, which boils
down to the case $X_{AB}^L(0)=0$ by the strong Markov property applied
at the first time $X_{AB}^L=0$.
\end{itemize}
Now, by Proposition~\ref{propeffectiverecomb} applied with $\phi_L$
replaced by $\sqrt{\phi_L}/2$, we have
\begin{eqnarray*}
\IP_{a_L'}\biggl[S^L \geq\frac{\psi_L}{2\sqrt{\phi_L}}
\biggr]&\leq&\IP
_{a_L'}\biggl[S^L \geq\frac{\sqrt{\phi_L}}{2}\biggl(1+\frac{\theta
\log
(\sqrt{\phi_L}\rho_L/2)}{r_L\rho_L}\biggr)\biggr]\\
& \leq& e^{-(\theta_1/2)\sqrt{\phi_L}}+e^{-((\theta-\theta
_2)/2)\sqrt
{\phi_L}\log(\sqrt{\phi_L}\rho_L/2)}.
\end{eqnarray*}
Substituting in (\ref{boundexcursion}) and using the asymptotic relation
$1-(1-e^{-t})^t\sim te^{-t}$ as $t\rightarrow\infty$, we obtain that
for $L$ large enough,
the quantity in $(\ref{effect1})$ is bounded by $\sqrt{\phi
_L}e^{-(\theta_1/4)\sqrt{\phi_L}}$.

As concerns $(\ref{effect2})$, observe that there exists $\theta_4>0$
such that for every
\mbox{$L\geq1$}, each of the $q_i^L-Q_i^L$ is stochastically bounded by an exponential
random variable with parameter $\theta_4$. Indeed, when $X_{AB}^L$ lies
within $B(0,(3/2)R_B)$,
the rate at which a coalescence occurs due to a large event is bounded
from below by a
positive constant. On the other hand, it is not difficult to check that
when $X_{AB}^L$ lies
within $B(0,(3/2)R_B)^c$, the rate at which the two lineages are sent
at a distance greater
than $3R_B$ by a large event is also bounded from below by a positive
constant. The quantity
in~(\ref{effect2}) is therefore bounded by
\[
\IP_{a_L}\Biggl[\sum_{i=1}^{\lfloor\sqrt{\phi_L}\rfloor
}(q_i^L-Q_i^L)\geq\frac{\psi_L}{2}\Biggr] \leq\IP\Biggl[\sum
_{i=1}^{\lfloor\sqrt{\phi_L}\rfloor}\mathcal{E}_i\geq\frac{\psi
_L}{2}\Biggr] \leq\exp\biggl\{-\frac{\psi_L}{2}+c\sqrt{\phi_L}
\biggr\},
\]
where $(\mathcal{E}_i)_{i\in\IN}$ is a sequence of i.i.d. exponential
random variables with
parameter~$\theta_4$ and $c$ is a positive constant expressed in terms
of the exponential
moment of~$\mathcal{E}_1$. The result follows.
\end{pf}

\subsection{Decorrelation time of two lineages starting in the same
individual}\label{subsdecorrelation}
In the previous section we obtained some information on the time
required for two lineages
starting in the same individual to become separated by a distance
greater than $3R_B$. We know that the lineages behave independently
whenever they are at distance greater than $2R_B$. However, nothing
guarantees that after the random time $\rmT_{(3R_B)}^L$ of Lemma~\ref
{lemmaseparation}, the ancestral lineages of $A$ and $B$ will evolve
independently. Indeed, it is very likely
that after some time they will once again be within distance $2R_B$ of
one another and
coalescence events will keep them close together for a potentially long
period of time.
Hence, in order to prove Theorem~\ref{theocorrelation}, we would like
to know how much time
our lineages need before they
start ``looking'' as if they were independent. That is, we are interested
in the time until their separation is of the same order as if they had
evolved according to
independent copies of $\ell^L$ started from $0$.
Recall from Lemma~\ref{lemA} that for (large) times less than $L^{2(1-\alpha
)}/\sqrt{\log L}$,
the difference of two independent lineages behaves like Brownian motion
on $\IR^2$.
The following proposition thus tells us that the decorrelation time we
are looking for is
asymptotically bounded from above by $(\log L)^5(1+\frac{\log
\rho
_L}{r_L\rho_L})$.
%
\begin{propn}\label{propdecorr}Let $(T_L)_{L\geq1}$ be a sequence of
times such that\break $ (\log L)^5(1+\frac{\log\rho_L}{r_L\rho
_L}
)\leq T_L\leq\frac{L^{2(1-\alpha)}}{\log L}$ for every $L$. Then,
\[
\lim_{L\rightarrow\infty}\IP_{a_L}\biggl[|X_{AB}^L(T_L)
|\notin
\biggl[\frac{\sqrt{T_L}}{\log L},\sqrt{T_L}\log L\biggr]\biggr] =0.
\]
\end{propn}

The scheme of the proof of Proposition~\ref{propdecorr} will again be
to decompose the
path of $X_{AB}^L$ into appropriate excursions and incursions. We shall
show that the
proportion of the time before $T_L$ that $X_{AB}^L$ spends in the
region of space where it
does \textit{not} evolve like the difference of two independent lineages
is asymptotically
negligible.

To this end, for every $L\in\IN$, let us define the stopping times
$(Q_i^L)_{i\geq0}$ and $(q_i^L)_{i\geq0}$ by $q_0^L=Q_0^L= 0$, and
for every $i\geq1$,
\begin{eqnarray*}
Q_i^L&:= &\inf\{t>q_{i-1}^L \dvtx  X_{AB}^L(t)\notin B(0,3R_B)\},
\\
q_i^L&:= &\inf\{t>Q_i^L \dvtx  X_{AB}^L(t)\in B(0,2R_B)\},
\end{eqnarray*}
with the convention that $\inf\varnothing=+\infty$. We also write $k_L$
for the number of
``excursions'' that start before time $T_L$, that is,
\[
k_L := \max \{i\dvtx  Q_i^L\leq T_L\}.
\]
The first step in proving Proposition~\ref{propdecorr} is to show that
%
\begin{lemma}\label{lemmanumberexc}
For every $\delta\in(0,1/2)$, there exist $K(\delta)>0$ such that for
all $L$ large enough,
\[
\IP_{a_L}[k_L > K(\delta)\log T_L]\leq\delta.
\]
\end{lemma}

We postpone the proof of Lemma~\ref{lemmanumberexc} until the end of
the section and instead
exploit it to prove Proposition~\ref{propdecorr}.

\begin{pf*}{Proof of Proposition~\protect\ref{propdecorr}}
We construct a coupling between $X_{AB}^L$ and a compound Poisson
process $Y^L$ which evolves
as the difference between two \textit{independent} copies of $\ell^L$.
Define $Y^L$ as follows: during an excursion of $X_{AB}^L$, $Y^L$ makes
the same jumps as $X_{AB}^L$ at the
same times, that is,
\[
\forall i\geq1, \forall t\in(Q_i^L,q_i^L],\qquad
Y^L(t)-Y^L(t-)=X_{AB}^L(t)-X_{AB}^L(t-).
\]
During the remaining time, $Y^L$ jumps independently of $X_{AB}^L$ with
a jump intensity
equal to twice that given in (\ref{intensities}) rescaled in an
appropriate manner.
It is easy to check that the law of $Y^L$ is indeed as claimed, since outside
$B(0,2R_B)$, $X_{AB}^L$ evolves like the difference of two independent
lineages and so
the jump intensity corresponding to the process $Y^L$ is equal to twice
that in the
rescaled version of~(\ref{intensities}) at any time. Furthermore, by
construction, the
difference between $X_{AB}^L$ and $Y^L$ changes only during the time intervals
$[q_{i-1}^L,Q_i^L]$. For convenience, we retain the notation $\IP$ for
the probability
measures on the (larger) space of definition of the pair
$(X_{AB}^L,Y^L)$, and
set $Y^L(0)=0$, $\IP_{a_L}$-a.s.

Let us call $I_L$ the amount of time before $T_L$ during which
$X_{AB}^L$ and $Y^L$ behave independently, that is,
\[
I_L:= \sum_{i=1}^{k_L}(Q_i^L-q_{i-1}^L)+(T_L-q_{k_L}^L)_+.
\]
If $\theta_2$ is as in Proposition~\ref{propeffectiverecomb}, we
have
%
%
\begin{eqnarray} \label{coupling1}
&&\IP_{a_L}\biggl[|X_{AB}^L(T_L)|\notin\biggl[\frac{\sqrt
{T_L}}{\log L},\sqrt{T_L}\log L\biggr]\biggr] \nonumber\\
& &\qquad\leq\IP_{a_L}\biggl[I_L< (\log L)^2 \biggl(1+\frac{2\theta_2
\log
(\rho_L\log L)}{r_L\rho_L}\biggr) ;
\nonumber
\\[-8pt]
\\[-8pt]
\nonumber
&&\hspace*{60pt}\qquad{} |X_{AB}^L(T_L)
|\notin\biggl[\frac{\sqrt{T_L}}{\log L},\sqrt{T_L} \log L
\biggr]
\biggr]\nonumber\\
& & \qquad\quad{} + \IP_{a_L}\biggl[I_L\geq(\log L)^2 \biggl(1+\frac{2\theta_2
\log( \rho_L\log L)}{r_L \rho_L}\biggr)\biggr].\nonumber
\end{eqnarray}
First, let us show that the second term in the right-hand side of (\ref
{coupling1})
converges to $0$ as $L\rightarrow\infty$. Let $\delta\in(0,1/2)$.
By Lemma~\ref{lemmanumberexc}, there exists $K>1$ such that for $L$
large enough,
$\IP_{a_L}[k_L>K\log T_L]\leq\delta$. Hence, we can write
\begin{eqnarray*}
&&\IP_{a_L}\biggl[I_L\geq(\log L)^2 \biggl(1+\frac{2\theta_2 \log
(\rho_L\log L)}{r_L \rho_L}\biggr)\biggr]\\
&&\qquad\leq \IP_{a_L}\biggl[I_L\geq(\log L)^2 \biggl(1+\frac{2\theta_2
\log
(\rho_L\log L)}{r_L \rho_L}\biggr); k_L\leq K\log T_L\biggr]
+\delta.
\end{eqnarray*}
Now, by the same reasoning as in (\ref{boundexcursion}), we have
%
%
\begin{eqnarray}\label{term1}
&&\IP_{a_L}\biggl[ I_L\geq(\log L)^2 \biggl(1+\frac{2\theta_2 \log
(
\rho_L\log L)}{r_L \rho_L}\biggr); k_L\leq K\log T_L
\biggr]\nonumber
\\
& &\qquad \leq\IP_{a_L}\Biggl[\sum_{i=1}^{\lfloor K\log T_L\rfloor
+1}Q_i^L-q_{i-1}^L\geq(\log L)^2\biggl(1+\frac{2\theta_2 \log(
\rho
_L \log L)}{r_L\rho_L}\biggr)\Biggr]
\nonumber
\\[-8pt]
\\[-8pt]
\nonumber
& &\qquad \leq1-\biggl(1-\sup_{a_L'} \IP_{a_L'}\biggl[Q_1^L\geq\frac
{(\log
L)^2}{K\log T_L +1}\\
&&\hspace*{120pt}\qquad{}\times\biggl(1+\frac{2\theta_2 \log(\rho_L\log
L
)}{r_L\rho_L}\biggr)\biggr]\biggr)^{\lfloor K \log T_L\rfloor
+1},\nonumber
\end{eqnarray}
where the supremum is taken over all initial configurations $a_L'$ in
which the distance
between the blocks containing $A$ and $B$ is at most $2R_B$. Again, as in
(\ref{boundexcursion}), we can restrict our attention to initial
configurations in which
$A$ and $B$ belong to the same block (recall from the proof of
Lemma~\ref{lemmaseparation} that
the rate at which a sequence of ``separating'' events occurs is bounded
from below by a positive
constant whenever $X_{AB}^L\neq0$). By assumption,
$\log T_L\leq2\log L$ and $K>1$, and so using Lemma~\ref{lemmaseparation} with
$\phi_L=(\log L)/(2K)$ for the last inequality we obtain that for all
large $L$,
uniformly in $a_L'$ as above,
\begin{eqnarray*}
&&\IP_{a_L'}\biggl[Q_1^L\geq\frac{(\log L)^2}{K\log T_L +1}
\biggl(1+\frac
{2\theta_2 \log(\rho_L\log L)}{r_L\rho_L}\biggr)\biggr]\\
& &\qquad\leq\IP_{a_L'}\biggl[Q_1^L\geq\frac{\log L}{2K}\biggl(1+\frac
{2\theta
_2 \log((2K)^{-1}\rho_L\log L)}{r_L\rho_L}\biggr)\biggr] \\
& & \qquad\leq\sqrt{(2K)^{-1}\log L} e^{-\theta_3\sqrt{(2K)^{-1}\log L}}.
\end{eqnarray*}
Consequently, we obtain from the asymptotic relation
$1-(1-te^{-t})^{t^2}\sim t^3e^{-t}$
that the quantity in the right-hand side of (\ref{term1}) tends to
zero as
$L\rightarrow\infty$ and
\[
\limsup_{L\rightarrow\infty} \IP_{a_L}\biggl[I_L\geq(\log L)^2
\biggl(1+\frac{2\theta_2 \log( \rho_L\log L)}{r_L \rho_L}
\biggr)
\biggr]\leq\delta.
\]
Since $\delta$ was arbitrary, this limit is actually zero.

Let us now show that the first term in the right-hand side of (\ref
{coupling1}) tends to zero as $L\rightarrow\infty$. To this end,
observe that it is bounded by
%
%
\begin{eqnarray} \label{term2}
&&\IP_{a_L}\biggl[I_L < (\log L)^2\biggl(1+\frac{2\theta_2 \log
(\rho
_L\log L)}{r_L\rho_L}\biggr); \nonumber\\
&&\qquad  |X_{AB}^L(T_L)-Y^L(T_L)|>(\log\log L)(\log
L)\biggl\{1+\frac{2\theta_2\log( \rho_L\log L)}{r_L\rho
_L}\biggr\}
^{1/2}\biggr]\\
&&\qquad{}+ \IP_{a_L}\biggl[I_L < (\log L)^2\biggl(1+\frac{2\theta_2 \log
(\rho_L\log L)}{r_L\rho_L}\biggr);
\nonumber
\\
&&\hspace*{35pt}\qquad |X_{AB}^L(T_L)
|\notin
\biggl[\frac{\sqrt{T_L}}{\log L},\sqrt{T_L}\log L\biggr];\nonumber\\
&&\hspace*{35pt}\qquad|X_{AB}^L(T_L)-Y^L(T_L)|\leq(\log\log L)(\log
L)\nonumber\\
&&\hspace*{142pt}\qquad{}\times\biggl\{1+\frac{2\theta_2\log(\rho_L\log L)}{r_L\rho
_L}\biggr\}
^{1/2}\biggr].\nonumber
\end{eqnarray}
Because the difference $X_{AB}^L-Y^L$ changes only during the periods
$[q_{i-1}^L,Q_i^L]$, during which $|X_{AB}^L|\leq3R_B$ and $Y^L$ jumps
around according to twice the jump intensity given by the appropriate
rescaling of (\ref{intensities}), the first term in (\ref{term2}) is
bounded by
\begin{eqnarray*}
&&\IP_{a_L}\biggl[I_L< (\log L)^2\biggl(1+\frac{2\theta_2 \log
(\rho
_L\log L)}{r_L\rho_L}\biggr);\\
& &\phantom{\IP_{a_L}\biggl[}   |\hat{Y}^L(I_L)|+3R_B>(\log\log L)(\log L)\biggl\{
1+\frac{2\theta_2 \log(\rho_L \log L)}{r_L\rho_L}\biggr\}
^{1/2}\biggr],
\end{eqnarray*}
where $\hat{Y}^L$ is an independent copy of $Y^L$ starting from $0$.
Hence, we also have as an upper bound
\[
\IP_{a_L}[|\hat{Y}^L(I_L)|>(\log\log L)\sqrt{I_L}-3R_B],
\]
which tends to zero by a standard use of Markov's inequality and
equation~(22) of~\cite{BEV2010}.

As concerns the second term in (\ref{term2}), it is bounded by
\begin{eqnarray*}
&&\IP_{a_L}\biggl[ |Y^L(T_L)|\notin\biggl[\frac{\sqrt
{T_L}}{\log
L}+(\log\log L)(\log L)\biggl\{1+\frac{2\theta_2\log(\rho
_L\log
L)}{r_L\rho_L}\biggr\}^{1/2}, \\
&&\hspace*{80pt}  \sqrt{T_L}\log L-(\log\log
L)(\log
L)\biggl\{1+\frac{2\theta_2\log(\rho_L\log L)}{r_L\rho
_L}\biggr\}
^{1/2}\biggr]\biggr] \\
& &\qquad =\IP_{a_L}\biggl[|Y^L(T_L)|\notin\biggl[\frac{\sqrt
{T_L}}{\log L} \bigl(1+\e_L^{(1)}\bigr),\sqrt{T_L}\log L \bigl(1-\e_L^{(2)}\bigr)
\biggr]\biggr],
\end{eqnarray*}
where by assumption on $T_L$ and the fact that $\rho_L\geq\log L$,
\[
\e_L^{(1)}:= \frac{(\log L)^2\log\log L}{\sqrt{T_L}}\biggl\{1+\frac
{2\theta_2 \log(\rho_L\log L)}{r_L\rho_L}\biggr\}
^{1/2}\leq C \frac{\log\log L}{\sqrt{\log L}}
\]
and
\[
\e_L^{(2)}:= \frac{\log\log L}{\sqrt{T_L}}\biggl\{1+\frac{2\theta
_2 \log
(\rho_L\log L)}{r_L\rho_L}\biggr\}^{1/2}\leq C'\frac{\log
\log
L}{(\log L)^{5/2}}.
\]
An application of the central limit theorem then gives the result.
\end{pf*}

The proof of Lemma~\ref{lemmanumberexc} rests upon the following lemma.
%
\begin{lemma}\label{lemmaexclength}There exists $C_q,v_q>0$ such
that for every $L$ large enough, $v_q\leq v \leq L^{2(1-\alpha)}/(\log
L)$ and every initial condition $a_L'$ in which the separation between
$A$ and $B$ belongs to $B(0,5R_B)\setminus B(0,3R_B)$,
\[
\IP_{a_L'}[q_1^L>v]\geq\frac{C_q}{\log v}.
\]
\end{lemma}

The proof of Lemma~\ref{lemmaexclength} uses the same arguments as
the second half
of the proof of Lemma~\ref{lemmaoutside0} (based on Skorokhod
embedding) and so we omit~it.

\begin{pf*}{Proof of Lemma~\protect\ref{lemmanumberexc}}
Our strategy is to show that if we choose $K$ large enough, the
probability that none of
the first $K \log T_L$ excursions outside $B(0,3R_B)$ has duration of
$\cO(T_L)$ is smaller
than $\delta$. To achieve this, let $K>0$. We have
\begin{eqnarray*}
\IP_{a_L}[k_L> K\log T_L]&=& \IP_{a_L}\bigl[Q^L_{\integ
{K\log
T_L}+1}\leq T_L\bigr] \\
&=&\IP_{a_L}\Biggl[\sum_{i=1}^{\integ{K\log T_L}}
(q_i^L-Q_i^L)+
\sum_{i=1}^{\integ{K\log T_L}+1}(Q_i^L-q_{i-1}^L) \leq
T_L
\Biggr] \\
&\leq& \IP_{a_L}\Biggl[\sum_{i=1}^{\integ{K\log T_L}}
(q_i^L-Q_i^L) \leq T_L\Biggr] \\
&\leq& \IP_{a_L}[\forall i\in\{1,\ldots,\integ{K\log T_L}\}
, q_i^L-Q_i^L\leq T_L].
\end{eqnarray*}
Using a recursion and Lemma~\ref{lemmaexclength} together with the
fact that $|X_{AB}^L(Q_i^L)|\in[3R_B,5R_B]$ (recall the jump lengths
are bounded by $2R_B$), we arrive at
\begin{eqnarray*}
\IP_{a_L}[\forall i\in\{1,\ldots,\integ{K\log T_L}\}, q_i^L-Q_i^L\leq T_L]&\leq& \biggl(1-\frac{C_q}{\log T_L}
\biggr)^{\integ{K\log T_L}}\\
& \rightarrow& e^{-KC_q} \qquad  \mbox{as } L\rightarrow
\infty.
\end{eqnarray*}
Now choose $K(\delta)$ large enough that $e^{-K(\delta)C_q}\leq
\delta/2$,
and Lemma~\ref{lemmanumberexc} is proved.
\end{pf*}

\subsection{Proof of the main results}\label{subsproofs}
Now that we understand decorrelation better, we can prove Theorems~\ref
{theocorrelation}
and~\ref{theosemi-correl}. Recall the rescalings of time by a factor~$\rho_L$ and of
space by $L^{-\alpha}$ that have been in force since the beginning of
Section~\ref{section2loci}
and the notation $\tau_{ij}^L$ for the coalescence time of lineages $i$
and $j$ in \textit{original}
units. In order to work in the rescaled setting, we define $\mathrm{t}
_{ij}^L:= \tau_{ij}^L/\rho_L$
for every $i,j\in\{A,a,B,b\}$, and $\mathrm{t}^L:=\mathrm
{t}_{Aa}^L\wedge\mathrm{t}_{Bb}^L$.
We denote the genealogical process (on the original space and time
scales) of the four loci
corresponding to step $L$ by~$\cA^L$. As explained in Section~\ref{results},
this Markov process takes its values in the set of all marked
partitions of $\{A,a,B,b\}$.
For any $t\geq0$, each block of $\cA^L(t)$ contains the labels of the
lineages present
in the same individual at (genealogical) time $t$, and its mark gives
the current
location on $\IT(L)$ of this common ancestor.

\begin{remark}\label{rempb}
Several times during the course of the proofs below we shall apply
Proposition~\ref{propdecorr} with $T_L=L^{2(\beta-\alpha)}$. Strictly
speaking,
we can only do this if
$L^{2(\beta-\alpha)}\geq(\log L)^5(1+\frac{\log\rho
_L}{r_L\rho
_L})$, at least\vadjust{\goodbreak}
for $L$ large enough, which is not guaranteed by (\ref{cond}). However,
if it is not the case,
we can still find a sequence $(\phi_L)_{L\in\IN}$ tending to infinity
and such that
\[
\phi_L L^{2(\beta-\alpha)}\geq(\log L)^5\biggl(1+\frac{\log\rho
_L}{r_L\rho_L}\biggr)\qquad  \forall L\in\IN\quad \mbox{and}\quad
\lim_{L\rightarrow\infty}\frac{\log\phi_L}{\log L}=0.
\]
Now, for the sake of clarity we presented the results of Lemma~\ref
{lemmameeting} at times of the form $\rho_LL^{2(t-\alpha)}$ but its
proof shows that, because $\log(\phi_LL^{2(\beta-\alpha)})\sim\log
(L^{2(\beta-\alpha)})$ as $L\rightarrow\infty$, we also have
\[
\lim_{L\rightarrow\infty}\IP_{a_L}\bigl[T_{Aa}^L>\rho_L\phi
_LL^{2(\beta
-\alpha)}\bigr]=1.
\]
(Another way to see this is to use the inequality
$\IP_{a_L}[T_{Aa}^L>\rho_L\phi_LL^{2(t-\alpha)}]\leq\IP
_{a_L}[T_{Aa}^L>\rho_L\phi_LL^{2(\beta-\alpha)}]$
for any fixed $t>\beta$ and $L$ large enough, and then let~$t$ tend to
$\beta$.)
Hence, all the above arguments carry over with $L^{2(\beta-\alpha)}$
replaced by
$\phi_LL^{2(\beta-\alpha)}$. Since the modifications are minor, we work
with $L^{2(\beta-\alpha)}$ in all cases.
\end{remark}

\begin{pf*}{Proof of Theorem~\protect\ref{theocorrelation}}
The main
difficulty is that we are
interested in the first coalescence times of the pairs $(A,a)$ and
$(B,b)$, regardless of
that of any other pair. As a consequence, several coalescence and subsequent
recombination events may occur before $\mathrm{t}^L$, creating some
correlation between
lineages originally far from each other ($A$ and $b$, e.g.).
The point is to show
that on the timescale of interest, decorrelation occurs fast enough for
the system of
ancestral lineages to behave like two independent genealogical
processes, one
for each locus.

Let us start by showing (a). Note that we can assume $\beta<1$,
since otherwise the result follows from Proposition~\ref{prop1locus}
and the bound
\begin{eqnarray}
\IP_{a_L}\bigl[\mathrm{t}^L\leq L^{2(1-\alpha)}\bigr]\leq\IP
_{a_L}\bigl[\mathrm{t}
_{Aa}^L\leq L^{2(1-\alpha)}\bigr]+ \IP_{a_L}\bigl[\mathrm
{t}_{Bb}^L\leq
L^{2(1-\alpha)}\bigr]\rightarrow0\nonumber\\
\eqntext{\mbox{as } L\rightarrow
\infty.}
\end{eqnarray}

Hence, suppose $\beta<1$, fix $t\in(\beta,1]$ (the case $t=\beta$ is
treated as above) and let $L\in\IN$. By the Markov property applied to
$\cA^L$ at time $\rho_LL^{2(\beta-\alpha)}$, we have
%
%
\begin{eqnarray}\label{bounda1}
&&\IP_{a_L}[\mathrm{t}^L>L^{2(t-\alpha)}]\nonumber\\
&&\qquad= \IE_{a_L}
\bigl[\ind_{\{\mathrm{t}
^L>L^{2(\beta-\alpha)}\}}\IP_{\cA^L(\rho_LL^{2(\beta-\alpha
)})}\bigl[\mathrm{t}
^L>L^{2(t-\alpha)}-L^{2(\beta-\alpha)}\bigr]\bigr]
\nonumber
\\[-8pt]
\\[-8pt]
\nonumber
&&\qquad=  \IE_{a_L}\bigl[\IP_{\cA^L(\rho_LL^{2(\beta-\alpha)})}
\bigl[\mathrm{t}
^L>L^{2(t-\alpha)}-L^{2(\beta-\alpha)}\bigr]\bigr]\\
& &\qquad\quad{}-\IE_{a_L}\bigl[\ind_{\{\mathrm{t}^L\leq L^{2(\beta-\alpha)}\}
}\IP_{\cA
^L(\rho_LL^{2(\beta-\alpha)})}\bigl[\mathrm{t}^L>L^{2(t-\alpha
)}-L^{2(\beta
-\alpha)}\bigr]\bigr].\nonumber
\end{eqnarray}
Again, the second term in (\ref{bounda1}) is bounded by
\[
\IP_{a_L}\bigl[\mathrm{t}_{Aa}^L\leq L^{2(\beta-\alpha)}\bigr] +
\IP_{a_L}
\bigl[\mathrm{t}_{Bb}^L\leq L^{2(\beta-\alpha)}\bigr],
\]
which tends to $0$ as $L\rightarrow\infty$ by Proposition~\ref{prop1locus}.
Since Lemma~\ref{lemmasep} shows that, with probability tending to
$1$, at most two lineages
at a time can meet at distance less than $2R_B$, we can define $\rmT
_1^L$ as the first
time two of the four lineages come within distance $2R_B$ of each other
and write
%
%
\begin{eqnarray}\label{bounda2}
&&\IE_{a_L}\bigl[\IP_{\cA^L(\rho_LL^{2(\beta-\alpha)})}
\bigl[\mathrm{t}
^L>L^{2(t-\alpha)}-L^{2(\beta-\alpha)}\bigr]\bigr]\nonumber\\
& &\qquad = \IE_{a_L}\bigl[\IP_{\cA^L(\rho_LL^{2(\beta-\alpha)})}
\bigl[\rmT
_1^L>L^{2(t-\alpha)}-L^{2(\beta-\alpha)}\bigr]\bigr]
\nonumber
\\[-8pt]
\\[-8pt]
\nonumber
& &\qquad\quad{}  + \IE_{a_L}\bigl[\IP_{\cA^L(\rho_LL^{2(\beta-\alpha
)})}
\bigl[\rmT_1^L\leq L^{2(t-\alpha)}-L^{2(\beta-\alpha)} ;\\
&&\hspace*{118pt}\qquad \mathrm{t}
^L>L^{2(t-\alpha)}-L^{2(\beta-\alpha)}\bigr]\bigr].\nonumber
\end{eqnarray}
Setting aside the first term in the right-hand side of (\ref{bounda2})
for a moment, we
further decompose the event corresponding to the second term:
%
%
\begin{eqnarray}\label{bounda3}
& &\IE_{a_L}\bigl[\IP_{\cA^L(\rho_LL^{2(\beta-\alpha)})}
\bigl[\rmT
_1^L\leq L^{2(t-\alpha)}-L^{2(\beta-\alpha)} ; \mathrm
{t}^L>L^{2(t-\alpha
)}-L^{2(\beta-\alpha)}\bigr]\bigr]\nonumber\\
& &\qquad =\IE_{a_L}\bigl[\IP_{\cA^L(\rho_LL^{2(\beta-\alpha)})}
\bigl[\rmT
_1^L\leq L^{2(t-\alpha)}-L^{2(\beta-\alpha)}; m_1^L\notin\{Aa,Bb\} ;\nonumber\\
&&\qquad\hspace*{180pt}
\mathrm{t}^L>L^{2(t-\alpha)}-L^{2(\beta-\alpha)}\bigr]\bigr]
\\
& & \qquad\quad{}   + \IE_{a_L}\bigl[\IP_{\cA^L(\rho_LL^{2(\beta-\alpha
)})}
\bigl[\rmT_1^L\leq L^{2(t-\alpha)}-L^{2(\beta-\alpha)} ; m_1^L\in\{
Aa,Bb\}
; \nonumber\\
&&\hspace*{190pt}\qquad\mathrm{t}^L>L^{2(t-\alpha)}-L^{2(\beta-\alpha)}
\bigr]
\bigr],\nonumber
\end{eqnarray}
where $m_1^L$ denotes the pair of labels of the lineages which ``meet''
at time $\rmT_1^L$.
Let us show that the second term in (\ref{bounda3}) tends to $0$ as
$L\rightarrow\infty$.
Using Lemma~\ref{lemmameeting}, we know that, with probability
tending to one, no pairs
of lineages starting at (rescaled) separation $L^{-\alpha}x_L$ have met
at distance less
than $2R_B$ by time $L^{2(\beta-\alpha)}$. Hence, until this time any
of these pairs taken
separately evolves like two independent compound Poisson processes, and
their mutual distance
at time $L^{2(\beta-\alpha)}$ lies within $[L^{\beta-\alpha}/(\log
L),L^{\beta-\alpha}\log L]$
with probability tending to one (by a standard application of the
Central Limit Theorem).
On the other hand, by condition~(\ref{cond}) we can use
Proposition~\ref{propdecorr}
with $T_L=L^{2(\beta-\alpha)}$ (see Remark~\ref{rempb}) and conclude
that with probability
tending to $1$, the distance at time $T_L$ between each pair of
lineages starting within
the same individual also lies in $[L^{\beta-\alpha}/(\log L),L^{\beta
-\alpha}\log L]$.
The situation has thus become rather symmetric by time $L^{2(\beta
-\alpha)}$.
Suppose, for instance, that $m_1^L=Aa$. Then, either $\rmT_1^L<
L^{2(t-\alpha)}-L^{2(\beta-\alpha)}-\log L$
and $\mathrm{t}^L>L^{2(t-\alpha)}-L^{2(\beta-\alpha)}$ or
$\rmT_1^L\in[L^{2(t-\alpha)}-L^{2(\beta-\alpha)}-\log
L,L^{2(t-\alpha
)}-L^{2(\beta-\alpha)}]$.
The probability of the first event tends to $0$ by (\ref{coaltime}),
which shows that
once $A$ and $a$ are gathered at distance smaller than $2R_B$, they
coalesce in a time smaller
than $\log L$.
Lemma~\ref{lemmahitting} (if $t<1$) or (\ref{hittingL2}) (if $t=1$)
shows that the
probability of the second event also tends to $0$ as $L\rightarrow
\infty$.
Hence, the second term in (\ref{bounda3}) does indeed vanish as
$L\rightarrow\infty$.

So far, we have obtained
%
%
\begin{eqnarray}\label{recursivearg}
&&\IP_{a_L}\bigl[\mathrm{t}^L>L^{2(t-\alpha)}\bigr]\nonumber\\
&&\qquad= \IE_{a_L}
\bigl[\IP_{\cA
^L(\rho_LL^{2(\beta-\alpha)})}\bigl[\rmT_1^L>L^{2(t-\alpha
)}-L^{2(\beta
-\alpha)}\bigr]\bigr]
\nonumber
\\[-8pt]
\\[-8pt]
\nonumber
& &\qquad\quad{} + \IE_{a_L}\bigl[\IP_{\cA^L(\rho_LL^{2(\beta-\alpha)})}
\bigl[\rmT
_1^L\leq L^{2(t-\alpha)}-L^{2(\beta-\alpha)} ; m_1^L\notin\{
Aa,Bb\} ;  \\
&&\hspace*{190pt}\qquad   \mathrm{t}^L>L^{2(t-\alpha
)}-L^{2(\beta
-\alpha)}\bigr]\bigr] + \delta_L^1,\nonumber
\end{eqnarray}
where $\delta_L^1\rightarrow0$ as $L\rightarrow\infty$. Next, by the
strong Markov property applied to $\cA^L$ at time $\rho_L\rmT_1^L$ and
the fact that $\rmT_1^L<\mathrm{t}^L$ a.s., we have
\begin{eqnarray*}
&&\IE_{a_L}\bigl[\IP_{\cA^L(\rho_LL^{2(\beta-\alpha)})}\bigl[\rmT
_1^L\leq
L^{2(t-\alpha)}-L^{2(\beta-\alpha)} ; m_1^L\notin\{Aa,Bb\} ; \\
&&\hspace*{145pt}\qquad\mathrm{t}
^L>L^{2(t-\alpha)}-L^{2(\beta-\alpha)}\bigr]\bigr] \\
&&\qquad =  \IE_{a_L}\bigl[\IE_{\cA^L(\rho_LL^{2(\beta-\alpha)})}
\bigl[\ind_{\{
\rmT_1^L\leq L^{2(t-\alpha)}-L^{2(\beta-\alpha)}; m_1^L\notin\{
Aa,Bb\}
\}}\\
&&\hspace*{97pt}\qquad{}\times \IP_{\cA^L(\rho
_L\rmT_1^L)}\bigl[\mathrm{t}^L>L^{2(t-\alpha)}-L^{2(\beta-\alpha
)}-\rmT_1^L
\bigr]\bigr]\bigr].
\end{eqnarray*}
If $t<1$, Lemma~\ref{lemmasep} tells us that with probability tending
to $1$, the mutual
distance between each of the $5$ pairs of lineages different from
$m_1^L$ at time $\rmT_1^L$
belongs to the interval $[(\rmT_1^L)^{1/2}/(\log L),(\rmT
_1^L)^{1/2}\log L]$.
If $t=1$, equation (\ref{hittingL2}) shows that we can replace
$\ind_{\{\rmT_1^L\leq L^{2(1-\alpha)}-L^{2(\beta-\alpha)}\}}$ by
$\ind_{\{\rmT_1^L\leq L^{2(1-\alpha)}/(\log L)\}}$, up to an
asymptotically vanishing error
term, and so Lemma~\ref{lemmasep} still applies. Hence, by the
uniform convergence stated
in Lemma~\ref{lemmameeting}, the probability that one of these pairs
meet at distance less
than $2R_B$ before $2\rmT_1^L$ tends to zero. Furthermore,
Proposition~\ref{propdecorr}
guarantees that with very high probability, the pair that meet at time
$\rmT_1^L$ is also at
a distance belonging to $[(\rmT_1^L)^{1/2}/(\log L),(\rmT
_1^L)^{1/2}\log L]$ after
another $\rmT_1^L$ units of time. (This statement uses a conditioning on
$\rmT_1^L$, which turns $2\rmT_1^L$ into a deterministic time and
enables us to use
Proposition~\ref{propdecorr}.) Defining $\rmT_2^L$ and $m_2^L$ in the
same manner as above
(we number the different quantities which appear here to make the
recursion clearer) and using
exactly the same arguments as those leading to~(\ref{recursivearg}),
we can thus write that
with probability tending to~$1$,
\begin{eqnarray*}
&&\IP_{\cA^L(\rho_L\rmT_1^L)}\bigl[\mathrm{t}^L> L^{2(t-\alpha
)}-L^{2(\beta
-\alpha)}-\rmT_1^L\bigr]\\
& &\qquad = \IE_{\cA^L(\rho_L\rmT_1^L)}\bigl[\IP_{\cA^L(\rho_L\rmT
_1^L)}
\bigl[\rmT_2^L> L^{2(t-\alpha)}-L^{2(\beta-\alpha)}-2\rmT_1^L
\bigr]\bigr]\\
& &\qquad\quad{} + \IE_{\cA^L(\rho_L\rmT_1^L)}\bigl[\IP_{\cA^L(\rho
_L\rmT
_1^L)}\bigl[\rmT_2^L\leq L^{2(t-\alpha)}-L^{2(\beta-\alpha)}-2\rmT
_1^L ;\\
&&\qquad\hspace*{124pt} m_2^L\notin\{Aa,Bb\} ;\\
&&\qquad\hspace*{128pt}   \mathrm{t}^L>
L^{2(t-\alpha
)}-L^{2(\beta-\alpha)}-2\rmT_1^L\bigr]\bigr] +\delta_L^2,
\end{eqnarray*}
with $\delta_L^2\rightarrow0$ as $L\rightarrow\infty$. It is easy to
check that the above equality is also valid if $L^{2(t-\alpha
)}-L^{2(\beta-\alpha)}-2\rmT_1^L\leq0$. By induction, we obtain for
any \mbox{$k\in\IN$}
%
%
\begin{eqnarray}\label{approxbyindependent}
& &\IP_{a_L}\bigl[\mathrm{t}^L>L^{2(t-\alpha)}\bigr]\nonumber\\
& &\qquad = \IE_{a_L}\bigl[\IP_{\cA^L(\rho_LL^{2(\beta-\alpha)})}
\bigl[\rmT
_1^L>L^{2(t-\alpha)}-L^{2(\beta-\alpha)}\bigr]\bigr] \nonumber\\
&&\qquad\quad{}+ \IE_{a_L}\bigl[\IE_{\cA^L(\rho_LL^{2(\beta-\alpha
)})}\bigl[\ind
_{\{\rmT_1^L\leq L^{2(t-\alpha)}-L^{2(\beta-\alpha)}; m_1^L\notin\{
Aa,Bb\}\}} \nonumber\\
&&\quad\qquad{}\times\IP
_{\cA^L(2\rho_L\rmT_1^L)}\bigl[\rmT_2^L> L^{2(t-\alpha
)}-L^{2(\beta
-\alpha)}-2\rmT_1^L\bigr]\bigr]\bigr]  + \cdots\nonumber\\
& &\qquad\quad{} + \IE_{a_L}\bigl[\IE_{\cA^L(\rho_LL^{2(\beta-\alpha)})}
\bigl[\ind
_{\{\rmT_1^L\leq L^{2(t-\alpha)}-L^{2(\beta-\alpha)};m_1^L\notin\{
Aa,Bb\}\}}\nonumber\\
&&\quad\qquad{}\times \IE_{\cA^L(2\rho_L\rmT_1^L)}\bigl[\ind_{\{\rmT_2^L
\leq L^{2(t-\alpha)}-L^{2(\beta-\alpha)}-2\rmT_1^L\}} \ind_{\{
m_2^L\notin\{Aa,Bb\}\}}\\
&&\qquad\quad{}\times \IE_{\cA^L(2\rho_L\rmT_2^L)}
\bigl[\cdots \IE_{\cA
^L(2\rho_L\rmT_{k-2}^L)}\bigl[\ind_{\{\rmT_{k-1}^L<L^{2(t-\alpha
)}-L^{2(\beta-\alpha)}-2\rmT_1^L-\cdots-2\rmT_{k-2}^L\}}\nonumber\\
&&\qquad\quad{}\times \ind_{\{m_{k-1}^L\notin\{Aa,Bb\}\}} \IP_{\cA
^L(2\rho_L\rmT_{k-1}^L)}\bigl[\rmT_k^L> L^{2(t-\alpha)}-L^{2(\beta
-\alpha
)}\nonumber\\
&&\hspace*{182pt}\qquad\quad{}-\cdots-2\rmT_{k-1}^L\bigr]\bigr]\cdots \bigr]\bigr]\bigr]
\bigr]\nonumber\\
&&\quad\qquad{} + \IE_{a_L}\bigl[\cdots \IP_{\cA^L(2\rho_L\rmT
_{k-1}^L)}\bigl[\rmT
_k^L\leq L^{2(t-\alpha)}-L^{2(\beta-\alpha)}-2\rmT_1^L-\cdots
-2\rmT
_{k-1}^L ;\nonumber\\
&&\hspace*{3pt}\qquad\quad{}\mathrm{t}^L>L^{2(t-\alpha
)}-L^{2(\beta-\alpha
)}-2\rmT_1^L-\cdots-2\rmT_{k-1}^L\bigr]\cdots \bigr]+\sum_{i=1}^k
\delta
_L^i,\nonumber
\end{eqnarray}
in which all occurrences of $L^{2(1-\alpha)}$ are replaced by
$L^{2(1-\alpha)}/(\log L)$ if we are considering the case $t=1$. In
order to stop the recursion, let us show that for any $\e>0$, there
exists $k\in\IN$ such that the last but one term in (\ref{approxbyindependent}) is bounded by $\e$ for all $L$ large enough. To this end,
define the sequence of random times $(\gamma_i^L)_{i\geq1}$ by
\[
\gamma_1^L:= \inf\bigl\{t\geq L^{2(\beta-\alpha)}\dvtx  2 \mbox
{ rescaled lineages meet at distance less than }2R_B\bigr\},
\]
and for any $i\geq2$,
\[
\gamma_{i}^L:= \inf\{t\geq2\gamma_{i-1}^L\dvtx  2 \mbox
{ rescaled lineages meet at distance less than }2R_B\}.
\]
A simple recursion shows that for all $i\in\IN$, $\gamma_i^L$ and
$2\gamma_i^L$ are stopping times. We can thus apply the strong Markov
property at time $\rho_L\gamma_1^L$, then $\rho_L\gamma_2^L$, and so
on, and obtain that
%
%
\begin{eqnarray}\label{comparison}
&&\IE_{a_L}\bigl[\IE_{\cA^L(\rho_LL^{2(\beta-\alpha)})}
\bigl[\ind_{\{\rmT
_1^L\leq L^{2(t-\alpha)}-L^{2(\beta-\alpha)}\}}\IE_{\cA^L(2\rho
_L\rmT
_1^L)}\bigl[\ind_{\{\rmT_2^L\leq L^{2(t-\alpha)}-L^{2(\beta-\alpha
)}-2\rmT_1^L\}} \nonumber\hspace*{-25pt}\\
& &\quad\qquad  \cdots \times\IP_{\cA^L(2\rho_L\rmT_{k-1}^L)}
\bigl[\rmT
_k^L\leq L^{2(t-\alpha)}-L^{2(\beta-\alpha)}-2\rmT_1^L-\cdots
-2\rmT
_{k-1}^L\bigr]\cdots\bigr]\bigr]\bigr]\hspace*{-25pt}\\
& &\qquad= \IP_{a_L}\bigl[\gamma_k^L\leq L^{2(t-\alpha)}\bigr].\nonumber\hspace*{-25pt}
\end{eqnarray}
Since with probability tending to $1$ at each time $2\gamma_i^L$ the
four lineages are at
distance of the order of $(\gamma_i^L)^{1/2}$ of each other,
Proposition~\ref{propkingman}
guarantees that, up to an asymptotically vanishing error term, the
conditional probability
that $\gamma_{i+1}^L$ is less than $L^{2(t-\alpha)}-L^{2(\beta
-\alpha
)}-2\gamma_1^L-\cdots-2\gamma_i^L$
is bounded from above by
$\mathcal{C}:=(1+c)(1-(\frac{\beta-\alpha}{t-\alpha}
)^6)$,
where $c>0$ can be chosen arbitrarily close to $0$. It remains to
choose $k\in\IN$ such that
$\mathcal{C}^k\leq\e$ and to notice that the left-hand side of~(\ref
{comparison}) is an
upper bound for the last but one term in (\ref{approxbyindependent})
to conclude.

Finally, let us show that the other terms in (\ref{approxbyindependent}) are close to those
corresponding to a system of four independent lineages. Using the
integer $k=k(\e)$ obtained
in the last paragraph, we rewrite the decomposition (\ref{approxbyindependent}) in terms
of $(\gamma_i^L)_{i\in\IN}$ as follows (we retain the notation $m_i^L$
for the labels of the two
lineages meeting at time $\gamma_i^L$ and we set $\gamma_0^L:=0$):
%
%
\begin{eqnarray}\label{decomp}
\IP_{a_L}\bigl[\mathrm{t}^L>L^{2(t-\alpha)}\bigr]&=& \eta_L(\e)+
\sum_{j=1}^k\IP
_{a_L}\bigl[\gamma_{j-1}^L\leq L^{2(t-\alpha)};\nonumber\\
&&\qquad\hspace*{54pt} m_l^L\notin\{Aa,Bb\}
 \ \forall l\in\{1,\ldots,j-1\};\\
& &  \hspace*{173pt}   \gamma_j^L>L^{2(t-\alpha)}
\bigr],\nonumber
\end{eqnarray}
where $\eta_L(\e)$ is the sum of the last but one term in (\ref{approxbyindependent}) and of
the error terms $\delta_L^i$, and is smaller than $2\e$ for $L$ large
enough by definition of
$k(\e)$. Now, let us denote by $\hat{\cA}^L$ a system of four
independent lineages moving around
on $\IT(L)$ according to the law of the motion of a single (unrescaled)
lineage, and let us
define $(\hat{\gamma}_i^L)_{i\geq1}$ in the same way as $(\gamma
_i^L)_{L\geq1}$ but with
$\cA^L$ replaced by $\hat{\cA}^L$. Let us also write $\hat{\mathrm
{t}}_{Aa}^L$
(resp., $\hat{\mathrm{t}}_{Bb}^L$) for the smallest time $t$ such
that the
lineages $A$ and $a$
(resp., $B$ and $b$) meet at distance less than $2R_BL^{\alpha}$ at
time $\rho_Lt$, and
$\hat{m}_i^L$ for the indices of the pair meeting at time $\hat
{\gamma}_i^L$.
Exactly the same chain of arguments as above leads to a decomposition of
$\IP_{a_L}[\hat{\mathrm{t}}_{Aa}^L\wedge\hat{\mathrm
{t}}_{Bb}^L>L^{2(t-\alpha)}]$
of the form
(\ref{decomp}), with another sequence $(\hat{\eta}_L(\e))_{L\geq1}$
whose terms are bounded
by $2\e$ whenever $L$ is large enough. Now, let us emphasize that
Proposition~\ref{propkingman} also applies to the meeting times at
distance less than
$2R_BL^{\alpha}$, before which the evolutions of $\cA^L$ and $\hat
{\cA
}^L$ have the same
distribution. As a consequence, morally, we should have that the
distributions of the pairs
of indices $m_i^L$ and $\hat{m}_i^L$ both converge to a uniform draw
from the set of distinct
pairs of labels (in other words, each pair has asymptotically the same
chance to be that
meeting), and, furthermore, if $\gamma_i^L$ and $\hat{\gamma}_i^L$ are
of the same logarithmic
order, so should $\gamma_{i+1}^L$ and $\hat{\gamma}_{i+1}^L$ be.

More formally, let us define, for every $L\in\IN$ and $j\geq1$,
\[
\cL_j^L:= \frac{\log\gamma_j^L}{2\log L} \ind_{\{\gamma_j^L\leq
L^{2(1-\alpha)}/(\log L)\}}+\infty \ind_{\{\gamma_j^L>
L^{2(1-\alpha
)}/(\log L)\}},
\]
and $\hat{\cL}_j^L$ in a similar manner. Our goal is to show that for
each $j$, the vectors $V_j^L:= (\cL_1^L,m_1^L,\ldots,\cL_j^L,m_j^L)$
and $\hat{V}_j^L:= (\hat{\cL}_1^L,\hat{m}_1^L,\ldots,\hat{\cL
}_j^L,\hat
{m}_j^L)$ converge in distribution as $L\rightarrow\infty$ to the same
random vector, whose law is obtained by successive uses of
Proposition~\ref{propkingman}. Thus, let us prove by recursion that
the distribution functions of the two vectors converge to the same
limit. The case $j=1$ is a direct consequence of Proposition~\ref{propkingman}, which shows that for any $s\in[\beta,1]$ and $i_1\neq i_2$,
\begin{eqnarray*}
\lim_{L\rightarrow\infty}\IP_{a_L}[\cL_1^L\leq s-\alpha ; m_1^L=i_1i_2] &=&\frac{1}{6}\biggl(1-\biggl(\frac{\beta-\alpha
}{s-\alpha}\biggr)^6\biggr)\quad \mbox{and}\\
\lim_{L\rightarrow\infty}\IP_{a_L}[\cL_1^L= \infty ; m_1^L=i_1i_2]&=& \frac{1}{6}\biggl(\frac{\beta-\alpha}{1-\alpha
}\biggr)^6.
\end{eqnarray*}
[Recall the analysis made at the beginning of the proof, according to
which the lineages meet before time $L^{2(\beta-\alpha)}$ with
probability tending to zero, and at that time they are all at pairwise
distance $\cO(L^{\beta-\alpha})$.]

Suppose the distribution functions of $V_j^L$ and $\hat{V}_j^L$
converge to the same (nondegenerate) limit as $L$ tends to infinity.
Let then $s\in[\beta,1]$, $i_1\neq i_2$ and $\mathcal{B}$ be an event
of the form $\{\cL_1^L\leq s_1-\alpha ; m_1^L=i_1^{(1)}i_2^{(1)} ; \ldots ; \cL_j^L\leq s_j-\alpha ; m_j^L=i_1^{(j)}i_2^{(j)}\}$ for
some given $\beta\leq s_1\leq\cdots\leq s_j\leq s$. Using the strong
Markov property with $\cA^L$ at time $2\rho_L\gamma_j^L=2\rho
_LL^{2\cL
_j^L}$ and recalling the definition of $T_1^L$ as the first time two
rescaled lineages come at distance less than $2R_B$ of each other, we
obtain
%
%
\begin{eqnarray}\label{convergenceV}
&&\IP_{a_L}[V_j^L\in\mathcal{B} ; \cL_{j+1}^L\leq s-\alpha ; m_{j+1}^L=i_1i_2] \nonumber\\
&&\qquad= \IE_{a_L}\Bigl[\ind_{\{V_j^L\in\mathcal{B}\}}\IP_{\cA^L(2\rho
_LL^{2\cL_j^L})}\bigl[T_1^L\leq L^{2(s-\alpha)}-2L^{2\cL_j^L}; m_1^L=i_1i_2\bigr]\Bigr]\nonumber\\
&&\qquad= \IE_{a_L}\biggl[\ind_{\{V_j^L\in\mathcal{B}\}}\times\frac
{1}{6}
\biggl(1-\biggl(\frac{\cL_j^L}{s-\alpha}\biggr)^6\biggr)\biggr] \\
& &\qquad\quad{}+ \IE_{a_L}\biggl[\ind_{\{V_j^L\in\mathcal{B}\}}\biggl\{\IP_{\cA
^L(2\rho
_LL^{2\cL_j^L})}[T_1^L\leq L^{2(s-\alpha)}-2L^{2\cL_j^L}; m_1^L=i_1i_2]\nonumber\\
& & \hspace*{197pt}\qquad\quad{}     - \frac{1}{6}\biggl(1-
\biggl(\frac{\cL
_j^L}{s-\alpha}\biggr)^6\biggr)\biggr\}\biggr].\nonumber
\end{eqnarray}
Since $V_j^L$ converges in distribution to $V_j^{\infty}$ as
$L\rightarrow\infty$, and since the law of $V_j^{\infty}$ does not
charge the boundary of $\mathcal{B}$, the first term in the right-hand
side of (\ref{convergenceV}) converges to
\[
\IE\biggl[\ind_{\{V_j^{\infty}\in\mathcal{B}\}}\times\frac
{1}{6}
\biggl(1-\biggl(\frac{\cL_j^{\infty}}{s-\alpha}\biggr)^6\biggr)\biggr]=: \IP
[V_j^{\infty}\in\mathcal{B}; \cL_{j+1}^{\infty}\leq s-\alpha; m_{j+1}^{\infty}=i_1i_2].
\]
For the second term in (\ref{convergenceV}), we already saw that, up
to an asymptotically
vanishing error term, we can insert the indicator function of the set
$\{\cA^L(2\rho_LL^{2\cL_j^L})\in\Gamma(L,4,\cL_j^L+\alpha)\}$ within
the expectation,
where $\Gamma(L,4,\eta)$ is defined at the end of Section~\ref{section1locus} as the set of all
configurations of four lineages in which all pairwise distances between
the locations of the
lineages belong to $[L^{\eta}/(\log L),L^{\eta}\log L]$. Now, we can
also replace the first
probability within the curly brackets by the probability that
$T_1^L\leq L^{2(s-\alpha)}$
and $m_1^L=i_1i_2$ by Lemma~\ref{lemmahitting}. Then, the uniform
convergence stated in
Proposition~\ref{propkingman} easily gives us that the second term in
the right-hand of
(\ref{convergenceV}) tends to $0$ as $L\rightarrow\infty$. Likewise,
as $L$ tends to infinity,
\begin{eqnarray*}
\IP_{a_L}[V_j^L\in\mathcal{B}; \cL_{j+1}^L=\infty; m_{j+1}^L=i_1i_2]&\rightarrow& \IE\biggl[\ind_{\{V_j^{\infty
}\in
\mathcal{B}\}} \frac{1}{6}\biggl(\frac{\cL_j^{\infty}}{1-\alpha
}
\biggr)^6\biggr]\\
&=: & \IP[V_j^{\infty}\in\mathcal{B}; \cL_{j+1}^{\infty
}=\infty; m_{j+1}^{\infty}=i_1i_2],
\end{eqnarray*}
and an analogous result can be established when we allow some of the
$\cL_i^L$, \mbox{$i\leq j$}
(and so the subsequent ones) to be infinite. Since this convergence
holds for all $s$ and
$i_1i_2$ as above, we obtain the convergence in law of $V_{j+1}^L$
toward $V_{j+1}^{\infty}$,
whose distribution is determined by the above limits. By the induction
principle, for every
$j\in\IN$ the sequence $(V_j^L)_{L\geq1}$ converges in distribution
to a random vector $V_j^{\infty}$.
Since the same arguments apply to $(\hat{V}_j^L)_{L\geq1}$, the
distribution function of
$\hat{V}_j^L$ also converges to that of $V_j^{\infty}$ and convergence
in distribution also holds.
As a consequence, coming back to (\ref{decomp}), we obtain that for
each term of the sum,
\begin{eqnarray*}
&& \bigl|\IP_{a_L}\bigl[\gamma_{j-1}^L\leq L^{2(t-\alpha)}; m_l^L\notin\{
Aa,Bb\}\ \forall l\in\{1,\ldots,j-1\}; \gamma_j^L>L^{2(t-\alpha
)}
\bigr] \\
& &\qquad{} -\IP_{a_L}\bigl[\hat{\gamma}_{j-1}^L\leq L^{2(t-\alpha)}; \hat
{m}_l^L\notin\{Aa,Bb\}\ \forall l\in\{1,\ldots,j-1\}; \\
&&\qquad\hspace*{234pt}\hat
{\gamma
}_j^L>L^{2(t-\alpha)}\bigr]\bigr|\rightarrow0
\end{eqnarray*}
as $L\rightarrow\infty$, and so
\[
\limsup_{L\rightarrow\infty}\bigl|\IP_{a_L}\bigl[\mathrm
{t}^L>L^{2(t-\alpha
)}\bigr]-\IP_{a_L}\bigl[\hat{\mathrm{t}}^L>L^{2(t-\alpha)}
\bigr]\bigr|\leq4\e.
\]
Since $\e$ was arbitrary, this limit is actually zero. But $\hat{\cA
}^L$ is a system of four independent lineages, and so
\[
\IP_{a_L}\bigl[\hat{\mathrm{t}}^L>L^{2(t-\alpha)}\bigr]=\IP
_{a_L}\bigl[\hat{\mathrm{t}
}_{Aa}^L>L^{2(t-\alpha)}\bigr]\times\IP_{a_L}\bigl[\hat{\mathrm{t}
}_{Bb}^L>L^{2(t-\alpha)}\bigr]\rightarrow\biggl(\frac{\beta-\alpha
}{t-\alpha}\biggr)^2
\]
by Proposition~\ref{prop1locus}. This concludes the proof of
Theorem~\ref{theocorrelation}(a).

The arguments for the case (b) are very similar, using this time
Lemma~\ref{lemB} for a bound on the probability that some lineages meet during a
small interval of time, Lemma~\ref{lemC} for the distance separating the other
lineages when two of them meet and merge and setting $\cL_j^L:= \gamma
_j^L/(\frac{1-\alpha}{2\pi\sigma^2}L^{2(1-\alpha)}\log L
)$.
\end{pf*}

The proof of Theorem~\ref{theosemi-correl} uses essentially the same
arguments, except that
now, before time $\rho_LL^{2(\gamma-\alpha)}$, we cannot use
Proposition~\ref{propdecorr} and the lineages starting within the
same individual are still highly correlated. In fact, because
recombination acts on a linear timescale whereas ancestral relations
evolve on an exponential timescale, the proof will show that a phase
transition occurs: during a first phase, recombination does not act and
so the ancestral lines of the two loci of the same individual are not
yet separated, and at time $\rho_LL^{2(\gamma-\alpha)}$ recombination
appears in the picture and is quick enough to fully decorrelate the
genealogies at the two loci.

\begin{pf*}{Proof of Theorem~\protect\ref{theosemi-correl}} The case
(a) is a consequence of the result for two lineages. Indeed, if
condition~(\ref{cond2}) is fulfilled, then necessarily $(\log\rho
_L)/(r_L\rho_L)$ tends to infinity and for any $\e>0$ there exists
$L_0(\e)$ such that for every $L\geq L_0(\e)$,
\[
\frac{\log\rho_L}{r_L\rho_L}\geq L^{2(\gamma-\alpha)-\e}.
\]
Hence, since we assumed $\rho_L\leq CL^{2\alpha}$, we have for $t\in
[\beta,\gamma)$, $\e:= \gamma-t$ and $L\geq L_0(\e)$,
\[
r_L\rho_L L^{2(t-\alpha)}\leq\log\rho_L  L^{2(t-\alpha-\gamma
+\alpha
)+(\gamma-t)}\leq C' L^{-(\gamma-t)}\log L\rightarrow0
\qquad\mbox{as } L\rightarrow\infty.
\]
Therefore, with probability tending to one, no recombinations occur by
time $\rho_LL^{2(t-\alpha)}$ and $\cA^L$ boils down to a system of two
lineages, one ancestral to each of the two individuals sampled.
Proposition~\ref{prop1locus} enables us to conclude.\looseness=1

If $t=\gamma$ and $r_L\rho_LL^{2(\gamma-\alpha)}$ does not tend to zero
(otherwise recombination is too slow and the same argument as above
applies), then the probability that there is no coalescence by time
$r_L^{-1}/(\log L)$ tends to $(\beta-\alpha)/(\gamma-\alpha)$. Indeed,
the recombination rate on the modified timescale is of the order of
$r_L\rho_L$, and so with high probability no recombinations separate
the two loci in any of our two sampled individuals before time
$(r_L\rho
_L)^{-1}/(\log L)$. Moreover,
\begin{eqnarray*}
\frac{\log((r_L\rho_L)^{-1}/(\log L))}{\log L}&=&\frac{\log
({\log\rho_L}/{(r_L\rho_L)})-\log\log\rho_L -\log\log
L}{\log
L}\\
&\rightarrow&2(\gamma-\alpha) \qquad \mbox{as } L\rightarrow
\infty,
\end{eqnarray*}
hence, by Proposition~\ref{prop1locus} (see also Remark~\ref{rempb}), the probability that no coalescence occurs before $r_L^{-1}/(\log
L)$ tends to $(\beta-\alpha)/(\gamma-\alpha)$. The last step is to
observe that, again by Proposition~\ref{prop1locus} and Remark \ref
{rempb}, the probability that any of the pairs of lineages $Aa$ and
$Bb$ (considered separately) coalesces during the time interval
$[r_L^{-1}/(\log L),\rho_LL^{2(\gamma-\alpha)}]$ tends to $0$ as $L$
tends to infinity.

For (b), apply the Markov property at time $\psi_L:= \rho_L
(L^{2(\gamma-\alpha)}\vee(\log L)^5(1+\frac{\log\rho
_L}{r_L\rho
_L}))$:
\begin{eqnarray*}
&&\IP_{a_L}\bigl[\tau_{Aa}^L\wedge\tau_{Bb}^L>\rho_LL^{2(t-\alpha
)}
\bigr]\\
&&\qquad= \IE_{a_L}\bigl[\ind_{\{\tau_{Aa}^L\wedge\tau_{Bb}^L>\psi_L\}
}\IP
_{\cA^L(\psi_L)}\bigl[\tau_{Aa}^L\wedge\tau_{Bb}^L>\rho
_LL^{2(t-\alpha
)}-\psi_L\bigr] \bigr] \\
&&\qquad = \frac{(\gamma-\alpha)^2}{(t-\alpha)^2} \IP_{a_L}[\tau
_{Aa}^L\wedge\tau_{Bb}^L>\psi_L] +o(1),
\end{eqnarray*}
where the second equality comes from Proposition~\ref{propdecorr},
Theorem~\ref{theocorrelation}(a) and dominated convergence. Now, by
the case (a) and Remark~\ref{rempb},
\[
\IP_{a_L}[\tau_{Aa}^L\wedge\tau_{Bb}^L>\psi_L
]\rightarrow\frac
{\beta-\alpha}{\gamma-\alpha}\qquad  \mbox{as } L\rightarrow
\infty,
\]
which yields the desired result.

Case (c) is identical to (b).
\end{pf*}

\section*{Acknowledgments}
We thank the referees for their very
careful reading and their useful remarks.



\printaddresses

\end{document}